\documentclass{elsart}

\usepackage{amssymb}
\usepackage{amsfonts}
\usepackage{amsmath}
\usepackage{mathrsfs}

\newcommand{\Levy} {L\'{e}vy~}

\newcommand{\ud}{\mathrm{d}}
\newcommand{\halmos}{\quad\hfill\mbox{$\Box$}}

\begin{document}
\begin{frontmatter}
\title{Conditions for certain ruin for the generalised Ornstein-Uhlenbeck process and the structure of the upper and lower bounds}
\author{Damien Bankovsky}
\ead{Damien.Bankovsky@maths.anu.edu.au}
\address{Mathematical Sciences Institute, Australian National University, Canberra, Australia}
\begin{abstract} For a bivariate \Levy process $(\xi_t,\eta_t)_{t\geq 0}$ the
generalised Ornstein-Uhlenbeck (GOU) process is defined as
\[V_t:=e^{\xi_t}\left(z+\int_0^t e^{-\xi_{s-}}\ud
\eta_s\right),~~t\ge0,\]where $z\in\mathbb{R}.$ We present
conditions on the characteristic triplet of $(\xi,\eta)$ which
ensure certain ruin for the GOU. We present a detailed analysis on
the structure of the upper and lower bounds and the sets of values
on which the GOU is almost surely increasing, or decreasing. This
paper is the sequel to \cite{BankovskySly08}, which stated
conditions for zero probability of ruin, and completes a significant
aspect of the study of the GOU.

\end{abstract}
\end{frontmatter}

\section{Notation and Theoretical Background}\label{intro section}

For a review of publications and applications for the GOU, see
\cite{BankovskySly08}. In Section 2 of this paper, we state results
on certain ruin for the GOU. Theorem 3.1 of Paulsen \cite{Paulsen98}
gives conditions for certain ruin for the GOU in the special case in
which $\xi$ and $\eta$ are independent. In \cite{BankovskySly08} it
is shown that this theorem does not hold for the general case.
Theorems \ref{paulsen theorem 1} and \ref{paulsen theorem 2} of
Section 2 give the required generalization, stated in terms of the
characteristic triplet of $(\xi,\eta).$ Section 3 begins with
results, in particular Proposition \ref{description of L
proposition} and Theorem \ref{combinations proposition}, which
describe the structure of the upper and lower bounds and the sets of
values on which the GOU is almost surely increasing, or decreasing.
Section 3 then outlines the ruin probability implications of these
structural results, in particular with Theorems \ref{minor certain
ruin theorem} and \ref{ruin prob theorem for L and U}, which state
conditions for certain ruin in terms of upper and lower bound
structure. Section 3 concludes with technical propositions used to
prove the major theorems. Section 4 contains proofs of the results
in Section 2 and 3, and concludes with a number of examples which
illustrate and extend certain results. For the remainder of this
section we set up some notation, which builds on that of
\cite{BankovskySly08}, and outline some basic results which we will
need.

Let $(\xi,\eta)$ be a bivariate \Levy process on a filtered complete
probability space $(\Omega, \mathscr{F},\mathbb{F},P)$ and define
the GOU process $V,$  and the associated stochastic integral process
$Z,$ as
\begin{equation}\label{GOU definition} V_t:=e^{\xi_t}\left(z+\int_0^t e^{-\xi_{s-}}\ud \eta_s\right),\end{equation}
and
\begin{equation}\label{Z definition} Z_t:=\int_0^t e^{-\xi_{s-}}\ud \eta_s.\end{equation} To avoid trivialities, assume that neither $\xi$ nor $\eta$ are
identically zero. It was shown in  \cite{BankovskySly08} that
\begin{equation}\label{jump equation 1}\Delta V_t=e^{\Delta\xi_t}\left(\Delta\eta_t-V_{t-}\left(e^{-\Delta\xi_t}-1\right)\right).
\end{equation}
The characteristic triplet of $(\xi,\eta)$ will be written
$\left((\tilde{\gamma}_\xi,\tilde{\gamma}_\eta),\Sigma_{\xi,\eta},\Pi_{\xi,\eta}\right).$
The characteristic triplet of $\xi$ as a one-dimensional \Levy
process will be written  $(\gamma_\xi,\sigma_\xi^2,\Pi_\xi),$ where
\begin{equation}\label{first 2 dim to 1 dim equation}\gamma_\xi=\tilde\gamma_\xi+
\int_{\{|x|< 1\}\cap\{x^2+y^2\ge 1\}}x\Pi_{\xi,\eta}(\ud(x,y)),
\end{equation}
and $\sigma_\xi^2$ is the upper left entry in the matrix $\Sigma_{\xi,\eta},$ and $\eta$ is symmetric.
The random jump measure and Brownian motion
components of $(\xi,\eta)$ will be denoted respectively by
$N_{\xi,\eta,t}$ and $(B_\xi,B_\eta).$

For a Lebesgue set $\Lambda$ define the \emph{hitting time} of
$\Lambda$ by $V$ to be
$T_{z,\Lambda}:=\inf\{t>0:V_t\in\Lambda|V_0=z\},$ where
$T_{z,\Lambda}:=\infty$ whenever $V_t\not\in\Lambda$ for all $t>0$
and $V_0=z.$ When the context makes it obvious we will simply write
$T_\Lambda$. Define the \emph{infinite horizon ruin probability} for
the GOU by
\[\psi(z):=P\left(\inf_{t>0}V_t<0|V_0=z\right)
=P\left(\inf_{t>0}Z_t<-z\right)
=P\left(T_{z,(-\infty,0)}<\infty\right).
\]
Note that for all $t>0,$ $V_t$ is increasing as a function of the initial value $z$ and hence, if $0\le z_1\le z_2,$ then $\psi(z_1)\ge\psi(z_2).$
For further explanation of the above terms, as well as extra definitions
and results for \Levy processes, see Section 1 of
\cite{BankovskySly08}. We now outline notation and theory needed for the present paper,
which were not dealt with in Section 1 of \cite{BankovskySly08}.

The \emph{total variation} of an $\mathbb{R}^n$-valued function over
the interval $[a,b]$ is defined by
\[V_f([a,b]):=\sup\sum_{i=1}^n\left|f\left(t_i\right)-f\left(t_{i-1}\right)\right|,
\] where the supremum is taken over all finite partitions
$a=t_0<t_1<\cdots<t_n=b.$ A \Levy process $X$ on $\mathbb{R}^n,$
with characteristic triplet $(\gamma_X,\Sigma_X,\Pi_X)$ and random
jump measure $N_{X,t},$ is said to be of \emph{finite variation} if,
with probability $1,$ its sample paths $X_t(\omega)$ are of finite
total variation on $[0,t]$ for every $t>0.$ It is shown in
\cite{ContTankov04}, p.86, this occurs iff $\Sigma_X=0$ and
$\int_{|z|\le 1}|z|\Pi_X(\ud z)<\infty.$ Further, if this occurs
then we can write
\[X_t=d_Xt+\int_{\mathbb{R}^n}zN_{X,t}(\cdot,\ud z)=dt+\sum_{0<s\le
t}\Delta X_s,
\] where
\begin{equation}\label{drift equation}d_X=\gamma_X-\int_{|z|<1}z\Pi_X(\ud z)\in\mathbb{R}^n=
E\left(X_1-\int_{\mathbb{R}^n}z N_{X,1}(\cdot,\ud
z)\right)\end{equation} is called the \emph{drift vector} of $X$. A
1-dimensional \Levy process $X$ is said to be a \emph{subordinator}
if $X_t(\omega)$ is an increasing function of $t,$ a.s., and it is
shown in \cite{ContTankov04}, p.88, that the following conditions
are equivalent:
\begin{enumerate}
\item $X$ is a subordinator.
\item $X_t\ge0$ a.s for some $t>0.$
\item $X_t\ge0$ a.s for every $t>0.$
\item The characteristic triplet satisfies
\[\sigma_X^2=0,~\int_{(-\infty,0]}\Pi_X(\ud
x)=0,~\int_{(0,1)}x\Pi_X(\ud x)<\infty,~\mathrm{and}~d_X\ge 0.
\] That is, there is no Brownian component, no negative jumps, the
positive jumps are of finite variation and the drift is
non-negative.
\end{enumerate}

A 1-dimensional \Levy process $X$ will drift to $\infty,$ drift to
$-\infty$ or oscillate between $\infty$ and $-\infty,$ namely, one of the following must hold:
\begin{equation}\label{drift to infinity}\lim_{t\rightarrow\infty}X_t=\infty~~\mathrm{a.s}.;
\end{equation}
\begin{equation}\label{drift to negative infinity}\lim_{t\rightarrow\infty}X_t=-\infty~~\mathrm{a.s.};
\end{equation}
\begin{equation}\label{oscillate}-\infty=\liminf_{t\rightarrow\infty}X_t<\limsup_{t\rightarrow\infty}X_t=\infty~~\mathrm{a.s.}
\end{equation}
Necessary and sufficient conditions for these cases are given in
\cite{DoneyMaller02}. Whenever the expected value of $X_1$ is a
well-defined member of the extended real numbers, cases (\ref{drift
to infinity}), (\ref{drift to negative infinity}), and
(\ref{oscillate}) equate respectively to $E(X_1)>0,$ $E(X_1)<0,$ and
$E(X_1)=0.$ For the case in which the expected value does not exist,
we need more notation. For $x>0,$ denote the tail functions of the
\Levy measure by
\[\overline{\Pi}_X^+(x):=\Pi_X((x,\infty)),~~\overline{\Pi}_X^-(x):=\Pi_X((-\infty,-x)),~~
\overline{\Pi}_X(x):=\overline{\Pi}_X^+(x)+\overline{\Pi}_X^-(x).\]
Define, for $x\ge 1,$
\[A_X^+(x):=\max\{\overline{\Pi}_X^+(1),1\}+\int_1^x\overline{\Pi}_X^+(u)\ud
u\] and
\[A_X^-(x):=\max\{\overline{\Pi}_X^-(1),1\}+\int_1^x\overline{\Pi}_X^-(u)\ud
u\]
 and define the integrals
\[J_X^+:=\int_1^\infty\left(\frac{x}{A_X^-(x)}\right)|\overline{\Pi}_X^+(\ud x)|~~\mathrm{and}~~
J_X^-:=\int_1^\infty\left(\frac{x}{A_X^+(x)}\right)|\overline{\Pi}_X^-(\ud x)|.\]
In \cite{DoneyMaller02} it is shown that if $E(X_1)$ is not well
defined, that is, if \[\int_1^\infty x\Pi_X(\ud
x)=\int_{-\infty}^{-1}|x|\Pi_X(\ud x)=\infty,\] then (\ref{drift to
infinity}) occurs iff $J_X^-<\infty,$ (\ref{drift to negative
infinity}) occurs iff $J_X^+<\infty$ and (\ref{oscillate}) occurs
iff $J_X^-=J_X^+=\infty.$

It is shown in \cite{CarmonaPetitYor01} that the GOU is a time
homogenous strong Markov process. In \cite{EricksonMaller05},
necessary and sufficient conditions are stated for a.s. convergence
of $Z_t$ to a finite random variable $Z_\infty$ as $t$ approaches
$\infty,$ whilst in \cite{LindnerMaller05}, necessary and sufficient
conditions are stated for stationarity of $V.$ We will need to use
these conditions, and to describe them we need some further
notation.

For a bivariate \Levy process $(X,Y)$ define the integral
\[I_{X,Y}:=\int_{(e,\infty)}\frac{\ln(y)}{A_X^+(\ln(y))}|\overline{\Pi}_Y(\ud y)|\]
and the auxiliary \Levy process $K^{X,Y}$ by
\[K^{X,Y}_t:=Y_t+\sum_{0<s\le t}\left(e^{\Delta X_s}-1\right)\Delta Y_s-t\mathrm{Cov}(B_{X,1},B_{Y,1}),\]
where Cov denotes the covariance. Theorem 2 of
\cite{EricksonMaller05} states that $Z_t$ converges a.s. to a finite
random variable $Z_\infty$ as $t\rightarrow\infty$ iff
$\lim_{t\rightarrow\infty}\xi_t=\infty$ a.s. and
$I_{\xi,\eta}<\infty.$ There is a special case in which, for some
$c\in\mathbb{R},$
\begin{equation}\label{degenerate case
equation}Z_t=c\left(e^{-\xi_t}-1\right)~~\mathrm{and}~~V_t=e^{\xi_t}(z-c)+c,
\end{equation} a.s. for all $t\ge 0.$
Exact conditions for this degenerate situation, given in terms of
the characteristic triplet of $(\xi,\eta),$ will be stated in
Proposition \ref{equivalence proposition}. In this situation,
$\lim_{t\rightarrow\infty}\xi_t=\infty$ a.s. implies that $Z_t$
converges a.s. to the constant random variable $Z_\infty=-c$ as
$t\rightarrow\infty,$ and in \cite{BertoinLindnerMaller07} it is
shown that this is the only case in which $Z_\infty$ is not a
continuous random variable. Note that, regardless of the asymptotic
behaviour of $\xi,$ if (\ref{degenerate case equation}) holds then
$V$ is strictly stationary iff $V_0=c.$ If (\ref{degenerate case
equation}) does not hold  for any $c\in\mathbb{R},$ then Theorem 2.1
of \cite{LindnerMaller05} states that $V$ is strictly stationary iff
the stochastic integral $\int_0^\infty e^{\xi_{s-}}\ud
K^{\xi,\eta}_s$ converges a.s. or, equivalently, iff
$\lim_{t\rightarrow\infty}\xi_t=-\infty$ a.s. and
$I_{-\xi,K^{\xi,\eta}}<\infty.$ In this case the stationary random
variable $V_\infty$ satisfies $V_\infty=_D\int_0^\infty
e^{\xi_{s-}}\ud K^{\xi,\eta}_s.$

\section{Conditions for Certain Ruin}\label{ruin probability section}

In Theorem 1 of \cite{BankovskySly08}, exact conditions were given
on the characteristic triplet of $(\xi,\eta)$ for the existence of
$u\ge 0$ such that $\psi(u)=0,$ and a precise value was given for
the value $\inf\{u\ge 0:\psi(u)=0\},$ where we use the convention
that $\inf\{\emptyset\cap[0,\infty)\}=\infty.$ It is a consequence
of Theorem \ref{paulsen theorem 1} below, that when the relevant
assumptions are satisfied, there exists $z\ge0$ such that
$\psi(z)<1$ iff there exists $u\ge 0$ such that $\psi(u)=0.$ Thus,
even though they are not stated explictly, Theorem \ref{paulsen
theorem 1} implies exact conditions on the characteristic triplet of
$(\xi,\eta)$ for certain ruin.

Statements (1) and (2) of Theorem \ref{paulsen theorem 1} are
generalizations to the dependent case of Paulsen's Theorem 3.1,
parts (a) and (b), respectively. Statement (1) of Theorem
\ref{paulsen theorem 1} also removes Paulsen's assumption of finite
mean for $\xi,$ and replaces his moment conditions with the precise
necessary and sufficient conditions for stationarity of $V.$ For
statement (2) of Theorem \ref{paulsen theorem 1}, a finite mean
assumption and moment conditions remain necessary.

\begin{thm}\label{paulsen theorem 1} Let $m:=\inf\{u\ge 0:\psi(u)=0\}.$
\begin{enumerate}
\item Suppose $\lim_{t\rightarrow\infty}\xi_t=-\infty$ a.s. and
$I_{-\xi,K^{\xi,\eta}}<\infty.$  Then $0<\psi(z)<1$ iff $0\le
z<m<\infty.$

\item Suppose $E(\xi_1)=0,$ $E(e^{|\xi_1|})<\infty$ and there exist $p,q>1$ with
$1/p+1/q=1$ such that $E\left(e^{-p\xi_1}\right)<\infty$ and
$E\left(|\eta_1|^q\right)<\infty.$ If, for all $c\in\mathbb{R},$ the
degenerate case (\ref{degenerate case equation}) does not hold, then
$0<\psi(z)<1$ iff $0\le z<m<\infty.$ If there exists
$c\in\mathbb{R}$ such that equation (\ref{degenerate case equation})
holds, then $\psi(z)<1$ iff $\psi(z)=0,$ which occurs iff $0\le c\le
z$.
\end{enumerate}
\end{thm}

\begin{rem}\label{paulsen theorem 1 remarks}
\begin{enumerate}
\rm{

\item In proving \cite{Paulsen98} Theorem 3.1 (b), Paulsen discretizes the GOU
at integer time points and then uses a recurrence result from
\cite{BabillotBougerolElie97}. His argument uses the inequality
$P(V_1<0|V_0=z)>0$ for all $z\ge 0,$ which is true in the
independent case if either $\xi$ or $\eta$ has a Brownian component,
or can have negative jumps. However, even in the independent case,
this inequality can fail to hold when $V_t$ decreases due to a
deterministic drift. For example, let $N$ and $M$ be independent
Poisson processes with parameter 1 and define $\xi_t:=-t+N_t$ and
$\eta_t:=-t+M_t.$ Let $T_z:=\inf\{t>0:V_t<0|V_0=z\}.$ Then $V_t\ge
(z+1)e^{-\xi_t}-1:=V_t'$ on $t\le T_z$ and $P(V_1'<0|V_0'=z)=0$
whenever $z>e^1-1.$  In proving statement (2) of Theorem
\ref{paulsen theorem 1} we get around this difficulty by
discretizing the GOU at random times $T_i$ and then showing that the
stated conditions result in $P(V_{T_1}<0|V_0=z)>0$ for all $z\ge 0$
in the general case.

\item Assume $\xi$ and $\eta$ are independent and $\eta$ is not a
subordinator. In this case, whenever $\xi$ drifts to $-\infty$ a.s.
or $\xi$ oscillates between $\infty$ and $-\infty$ a.s., it is a
consequence of Theorem 1 in \cite{BankovskySly08}, that $\psi(u)>0$
for all $u\ge 0,$ and hence $m=\infty.$ Thus, by statement (1) of
Theorem \ref{paulsen theorem 1}, if
$\lim_{t\rightarrow\infty}\xi_t=-\infty$ a.s. and
$I_{-\xi,K^{\xi,\eta}}<\infty,$ then $\psi(z)=1$ for all $z\ge 0.$
This result is a slight strengthening of Paulsen's Theorem 3.1 (a).
Further, statement (2) simplifies exactly to Paulsen's Theorem 3.1
(b). Since $\xi$ and $\eta$ are independent the conditions in
statement (2) simplify to $E(\xi_1)=0,$
$E\left(e^{|\xi_1|}\right)<\infty$ and $E(\eta_1)<\infty.$ Since
$m=\infty,$ $\psi(z)=1$ for all $z\ge0$ whenever these conditions
hold. The simplification of conditions occurs because H\"{o}lder's
inequality is not needed in the proof, and a simpler argument using
independence suffices. When transferred onto the \Levy measure,
these conditions are equivalent to those in Paulsen's Theorem 3.1
(b). }
\end{enumerate}
\end{rem}

We now present Theorem \ref{paulsen theorem 2}, which is the
generalization to the dependent case of Paulsen's Theorem 3.1, part
(c). In addition, Paulsen's assumption of finite mean for $\xi$ is
removed, and his moment conditions are replaced with the precise
necessary and sufficient conditions for a.s. convergence of $Z_t$ to
a finite random variable $Z_\infty,$ as $t\rightarrow\infty.$ A
formula for the ruin probability in this situation was given in
Theorem 4 of \cite{BankovskySly08}, however no conditions for
certain ruin were found. Theorem \ref{paulsen theorem 2} gives exact
conditions on the characteristic triplet of $(\xi,\eta)$ for certain
ruin. To state these conditions, we need the following definitions.

Let $A_1:=\left\{(x,y)\in\mathbb{R}^2:x\ge 0,y\ge0\right\},$ and
similarly, let $A_2,$ $A_3$ and $A_4$ be the quadrants in which
$\{x\ge0,y\le0\},$ $\{x\le0,y\le0\}$ and $\{x\le0,y\ge0\}$
respectively. For each $i=1,2,3,4$ and $u\in\mathbb{R}$ let
\[B_i^u:=\left\{(x,y)\in A_i:y-u(e^{-x}-1)>0\right\}\] and define
\[
\theta_1':= \left\{ \begin{array}{ll}
\inf\left\{u\le 0:\Pi_{\xi,\eta}(B_1^u)>0\right\} & \\
0~~~~\textrm{if~}\Pi_{\xi,\eta}(A_1\setminus A_2)=0, &
\end{array} \right.\,
\theta_3':= \left\{ \begin{array}{ll}
\sup\left\{u\le 0:\Pi_{\xi,\eta}(B_3^u)>0\right\} &\\
-\infty~~~~\textrm{if~}\Pi_{\xi,\eta}(A_3\setminus A_2)=0,&
\end{array} \right.\,
\]
\[
\theta_2':= \left\{ \begin{array}{ll}
\inf\left\{u\ge 0:\Pi_{\xi,\eta}(B_2^u)>0\right\} & \\
\infty~~~~\textrm{if~}\Pi_{\xi,\eta}(A_2\setminus A_3)=0, &
\end{array} \right.\,
\theta_4':= \left\{ \begin{array}{ll}
\sup\left\{u\ge 0:\Pi_{\xi,\eta}(B_4^u)>0\right\} &\\
0~~~~\textrm{if~}\Pi_{\xi,\eta}(A_4\setminus A_3)=0. &
\end{array} \right.\,
\]

\begin{thm}\label{paulsen theorem 2} Suppose $\lim_{t\rightarrow\infty}\xi_t=\infty$ a.s. and
$I_{\xi,\eta}<\infty.$ Then $\psi(0)=1$ if and only iff $-\eta$ is a
subordinator, or there exists $z>0$ such that $\psi(z)=1.$ The
latter occurs if and only if $\Pi_{\xi,\eta}(A_1)=0,$
$\theta_4'\le\theta_2',$ and there exists
$u\in[\theta_4',\theta_2']$ such that
\begin{equation}
\label{covariance matrix equation} \Sigma_{\xi,\eta}= \left[
  \begin{array}{ r r }
     1 & -u \\
     -u & u^2
  \end{array} \right]\sigma_\xi^2,
  \end{equation} and
  \begin{equation}\label{finite drift equation}g(u):=\tilde{\gamma}_\eta+u\tilde{\gamma_\xi}-\frac{1}{2}u\sigma_\xi^2-
\int_{\{x^2+y^2<1\}}(ux+y)\Pi_{\xi,\eta}(\ud(x,y))\le0.
\end{equation}
If there exists $z\ge0$ such that $\psi(z)=1$ and, for all
$c\in\mathbb{R},$ the equation (\ref{degenerate case equation}) does
not hold, then the following hold:
\begin{enumerate}
\item If $\sigma_\xi^2=0$ then $\psi(z)=1$ for all $z\le
m:=\sup\left\{u\in[\theta_4',\theta_2']:g(u)\le 0 \right\},$ and
$0\le\psi(z)<1$ for all $z>m;$

\item If $\sigma_\xi^2\neq0$ then $\psi(z)=1$
for all $z\le m:=-\frac{\sigma_{\xi,\eta}}{\sigma_\xi^2},$ and
$0<\psi(z)<1$ for all $z>m.$
\end{enumerate}

If there exists $z\ge0$ such that $\psi(z)=1$ and there exists
$c\in\mathbb{R}$ such that (\ref{degenerate case equation}) holds,
then $0<c=\theta_4'=\theta_2',$ $\psi(z)=1$ for all $z<c,$ and
$\psi(z)=0$  for all $z\ge c.$
\end{thm}

\begin{rem}\label{paulsen theorem 2 remarks}
\begin{enumerate}
\rm{
\item When $\Pi_{\xi,\eta}(A_1)=0,$ $\theta_4'\le\theta_2'$ and
$u\in[\theta_4',\theta_2']$ the function $g(u)$ is a well-defined
member of the extended reals. The existence and finiteness of $g$ is
fully analysed in point (1) of Remark \ref{lower bound remarks}.

\item Assume $\xi$ and $\eta$ are independent. Then all jumps occur at the axes
of the sets $A_i,$ and $\sigma_{\xi,\eta}=0.$ With a little work,
Theorem \ref{paulsen theorem 2} simplifies to the following
statement: Suppose $\lim_{t\rightarrow\infty}\xi_t=\infty$ a.s. and
$I_{\xi,\eta}<\infty.$ Then $\psi(0)=1$ iff $-\eta$ is a
subordinator, or $\psi(z)=1$ for some $z>0.$ The latter occurs iff
$\xi$ and $\eta$ are each of finite variation and have no positive
jumps, and $g(z)\le 0.$ Note that when $(\xi,\eta)$ is finite
variation, $g$ simplifies to $g(u)=d_\eta+ud_\xi,$ as explained in
equation (\ref{linear g equation}). Since $\xi$ drifts to $\infty$
a.s., it must be that $d_\xi>0.$ Thus, $g(z)\le 0$ for some $z>0$
iff $d_\eta<0.$ In particular, $-\eta$ is a subordinator.

\item In Paulsen \cite{Paulsen98}, Theorem 3.1 (c), it is stated that when $\xi$ and
$\eta$ are independent, $E(\xi_1)>0,$ and a set of moment conditions
hold, then $\psi(z)=1$ iff $\xi_t=\alpha t,$ $\eta_t=\beta t$ and
$\beta<-\alpha z$ for real constants $\alpha$ and $\beta.$ This
statement contradicts the independence version of Theorem
\ref{paulsen theorem 2} stated above, and is false. A simple
counterexample is $(\xi,\eta)_t:=(t,-t-N_t)$ where $N$ is a Poisson
process. Paulsen's moment conditions are satisfied trivially.
However, Theorem \ref{paulsen theorem 2} implies that $\psi(z)=1$
for all $z\le 1,$ and this is confirmed by elementary calculations.
If we denote the jump times of $N_t$ by $0=T_0<T_1<T_2<\cdots$ then
\[V_t=1+e^t\left(z-1-\sum_{i=1}^{N_t}e^{-T_i}\right).
\] Thus, if $z=1,$ then $V_{T_2}=-e^{T_2-T_1}<0$ a.s. and so $\psi(1)=1.$

}
\end{enumerate}
\end{rem}

The following proposition fully explains the ruin probability
function for the degenerate situation (\ref{degenerate case
equation}). It will be used to prove that Theorems \ref{paulsen
theorem 1} and \ref{paulsen theorem 2} correctly allow for this
case.

\begin{prop}\label{degenerate certain
ruin prop} Suppose that there exists $c\in\mathbb{R}$ such that
$V_t=e^{\xi_t}(z-c)+c.$ If $c\ge 0$ then $\psi(z)=0$ for all $z\ge
c,$ and the following statements hold for all $0\le z<c:$
\begin{enumerate}
\item If $\xi$ drifts to $-\infty$ a.s. then $0<\psi(z)<1;$
\item If $\xi$ oscillates between $\infty$ and $-\infty$ a.s. then $\psi(z)=1;$
\item If $\xi$ drifts to $\infty$ a.s. then $\psi(z)=1.$
\end{enumerate}
If $c<0$ then the following statements hold for all $z\ge 0:$
\begin{enumerate}
\item[(4)] If $\xi$ drifts to $-\infty$ a.s. then $\psi(z)=1;$
\item[(5)] If $\xi$ oscillates between $\infty$ and $-\infty$ a.s. then $\psi(z)=1;$
\item[(6)]  If $\xi$ drifts to $\infty$ a.s. then $0<\psi(z)<1.$
\end{enumerate}
\end{prop}

\section{Structure of the upper and lower bounds, and relationship with certain ruin}\label{bounds
section}

Define the \emph{lower bound function} $\delta$ and the \emph{upper bound function} $\Upsilon$ by
\[
\delta(z):= \inf \left\{u\in\mathbb{R}:P\left(\inf_{t\geq 0} V_t
\leq u \big|V_0=z\right) > 0 \right\}
\]
and
\[
\Upsilon(z):=\sup \left\{u\in\mathbb{R}:P\left(\sup_{t\geq 0} V_t
\geq u \big|V_0=z\right) > 0 \right\},
\] where we use the convention that $\inf \{\emptyset\cap\mathbb{R}\}=\infty$ and $\sup\{\emptyset\cap\mathbb{R}\}=-\infty.$
When $V_0=z,$ the probability that the sample paths $V_t$ will ever
rise above $\Upsilon(z),$ or below $\delta(z),$ is zero. In
particular, the ruin probability function $\psi$ satisfies
$\psi(z)=0$ iff $\delta(z)\ge0.$ Define the sets $L$ and $U$ by
\[L:=\{u\in\mathbb{R}:\delta(u)=u\}~~\mathrm{and}~~U:=\{u\in\mathbb{R}:\Upsilon(u)=u\}.
\]
It will be a consequence of Proposition \ref{endpoints of L
proposition} that $L$ and $U$ must each be of the form
\begin{equation}\label{types of L and U equation}\emptyset,\{a\},[a,b],[a,\infty),~\mathrm{or~}(-\infty,b]
\end{equation} for some $a,b\in\mathbb{R}.$ The fact that $L$ and $U$ are both connected sets is of great importance.

This section contains a detailed analysis of $\delta,$ $\Upsilon,$
$U$ and $L$ and their relationship with the ruin function. In
particular, we are interested in which combinations of $L$ and $U$
can exist. For each combination we are also interested in the
possible asymptotic behaviour of $\xi,$ namely, whether $\xi$ drifts
to $\infty$ a.s., $\xi$ drifts to $-\infty$ a.s. or $\xi$ oscillates
between $\infty$ and $-\infty$ a.s. We are interested in this
asymptotic behaviour because of its link with the conditions for
convergence of $Z_t$ and stationarity of $V,$ as discussed in
Section \ref{intro section}. As well as being of independent
interest, the results contained in this section are essential for
the proofs of Theorems \ref{paulsen theorem 1} and \ref{paulsen
theorem 2}.

We begin with comments on $\delta,$ and $L.$ The analogues for
$\Upsilon$ and $U$ are obvious through symmetry. Firstly, note that
$\delta(z)\le z$ for all $z\in\mathbb{R},$ whilst the fact that
$V_t$ is increasing in $z$ for all $t\ge 0$ implies that
$\delta(z_1)\le\delta(z_2)$ whenever $z_1<z_2.$ The following
proposition explains the behaviour of the $\delta$ outside the set
$L,$ and states that $L$ is precisely the set of starting parts
$V_0=z$ for which almost all sample paths $V_t$ are increasing for
some time period. Recall that
$T_{z,\Lambda}:=\inf\{t>0:V_t\in\Lambda\},$ and define
$L^c:=\mathbb{R}\setminus L.$

\begin{prop}\label{description of L proposition}The following statements hold for $L$ and $\delta$,
and the symmetric statements hold for $U$ and $\Upsilon$:
\begin{enumerate}
\item If $z\ge\sup L$ then $\delta(z)=\sup L;$
\item If $z<\inf L$ then $\delta(z)=-\infty;$
\item For $z\in L,$ $P\left(V_t\mathrm{~is~increasing~on~}0<t\le T_{z,L^c}|~V_0=z\right)=1;$
\item For $z\in L^c,$ $P\left(V_t\mathrm{~is~increasing~on~}0<t\le T_{z,L}|~V_0=z\right)<1.$
\end{enumerate}
\end{prop}

In Section \ref{intro section} we assumed that neither $\xi$ nor
$\eta$ are identically zero in order to avoid trivialities. The
following proposition explains these trivialities.
\begin{prop}\label{degenerate proposition}
\begin{enumerate}
\item $L=\mathbb{R}$ iff $\xi_t=0$ a.s. for all $t>0$ and $\eta$ is a subordinator.
\item $U=\mathbb{R}$ iff $\xi_t=0$ a.s. for all $t>0$ and $-\eta$ is a subordinator.
\item $L=U=\mathbb{R}$ iff $\xi_t=\eta_t=0$ a.s. for all $t>0.$
\end{enumerate}
\end{prop}
For the rest of this paper we again assume that neither $\xi$ nor
$\eta$ are identically zero. The following proposition explains the
degenerate situation described in equation (\ref{degenerate case
equation}). Note that the deterministic case
$(\xi,\eta)_t:=(\alpha,\beta)t$ for non-zero constants $\alpha$ and
$\beta$ satisfies the conditions of this proposition for
$c=-\beta/\alpha.$ Recall that a Borel set
$\Lambda\subsetneq\mathbb{R}$ is an \emph{absorbing} set for $V,$ if
for all $0\le s\le t,$ $P(V_t\in\Lambda|V_s=x)=1$ for all
$x\in\Lambda.$ That is, whenever a sample path $V_t$ hits $\Lambda,$
it never leaves. The stochastic exponential will be denoted by
$\epsilon.$
\begin{prop}\label{equivalence proposition} The following are equivalent for $c\neq 0$:
\begin{enumerate}
\item $L\cap U\neq\emptyset$;
\item $L\cap U=\{c\}$;
\item $V_t=e^{\xi_t}(z-c)+c$ and $Z_t=c\left(e^{-\xi_t}-1\right);$
\item $\{c\}$ is an absorbing set;
\item $\Sigma_{\xi,\eta}$ satisfies (\ref{covariance matrix
equation}) for $u=c$, $\Pi_{\xi,\eta}=0$ or is supported on the
curve $\{(x,y):y-c(e^{-x}-1)=0\},$ and $g(c)=0;$
\item $e^{-\xi_t}=\epsilon(\eta/c)_t.$
\end{enumerate}
If the above conditions hold and $\Sigma_{\xi,\eta}\neq0$ then
$L=U=\{c\}$ and there exist \Levy processes $(\xi,\eta)$ for this
situation such that $\xi$ drifts to $\infty$ a.s., $\xi$ drifts to
$-\infty$ a.s. or $\xi$ oscillates a.s.  If the above conditions
hold and $\Sigma_{\xi,\eta}=0$ then:
\begin{enumerate}
\item[(a)] $U=(-\infty,c]$ and $L=[c,\infty)$ iff $\xi$ is a
subordinator;
\item[(b)] $L=(-\infty,c]$ and $U=[c,\infty)$ iff $-\xi$ is a
subordinator;
\item[(c)] $L=U=\{c\}$ iff neither $\xi$ or $-\xi$ is a
subordinator. There exist \Levy processes $(\xi,\eta)$ for this
situation such that $\xi$ drifts to $\infty$ a.s., $\xi$ drifts to
$-\infty$ a.s. or $\xi$ oscillates a.s.
\end{enumerate}
\end{prop}

We present a theorem which describes all possible combinations of
$L$ and $U$ and the associated asymptotic behaviour of $\xi,$ for
the case in which $L\cap U=\emptyset.$

\begin{thm}\label{combinations proposition} Suppose that $L\cap U=\emptyset.$ If $\Sigma_{\xi,\eta}\neq0$
then only the following cases can exist:
\begin{enumerate}
\item $L=U=\emptyset$;
\item $L=\{a\}$ for some $a\in\mathbb{R}$ and $U=\emptyset$;
\item $U=\{a\}$ for some $a\in\mathbb{R}$ and $L=\emptyset$.
\end{enumerate}

If $\Sigma_{\xi,\eta}=0$ then only the following cases can exist:
\begin{enumerate}
\item[(a)] If $L=\emptyset$ then $U$ is of the form
$\emptyset,$ $\{a\},$ $[a,b],$ $[a,\infty),$ or $(-\infty,b]$ for
some $a,b\in\mathbb{R};$
\item[(b)] If $U=\emptyset$ then $L$ is of the form
$\emptyset,$ $\{a\},$ $[a,b],$ $[a,\infty),$ or $(-\infty,b]$ for
some $a,b\in\mathbb{R};$
\item[(c)] If $L\neq\emptyset$ and $U\neq\emptyset$ then there exist $a<b$ such that
$L=(-\infty,a]$ and $U=[b,\infty),$ or $U=(-\infty,a]$ and
$L=[b,\infty).$
\end{enumerate}
If $U=(-\infty,a]$ or $L=[b,\infty)$ (or both) then $\xi$ is a
subordinator. If $L=(-\infty,a]$ or $U=[b,\infty)$ (or both) then
$-\xi$ is a subordinator. For all of the other combinations of $L$
and $U$ above, there exist \Levy processes $(\xi,\eta)$ such that
$\xi$ drifts to $\infty$ a.s., $\xi$ drifts to $-\infty$ a.s. or
$\xi$ oscillates a.s.
\end{thm}

An absorbent set $\Lambda\subsetneq\mathbb{R}$ is a \emph{maximal
absorbing} set if it is not properly contained in any other
absorbing set. Note that if $\Lambda$ is a maximal absorbing set,
then $\mathbb{R}\setminus\Lambda$ contains no absorbing sets
otherwise we could take the union of $\Lambda$ with the absorbing
set, and this would be an absorbing set properly containing
$\Lambda.$ The following corollary is immediate. For each statement
(1)-(4), the claim that the sets $\Lambda$ are maximal absorbing
follows from Proposition \ref{description of L proposition}. The
remaining statements follow immediately from Theorem
\ref{combinations proposition}.

\begin{cor}\label{max absorbent proposition}
There exist \Levy processes $(\xi,\eta)$ with $L\cap U=\emptyset$
such that the associated GOU has the following maximal absorbing
sets $\Lambda:$
\begin{enumerate}
\item $\Lambda=U\cup L,$ where $U=(-\infty,a]$ and
$L=[b,\infty)$;
\item $\Lambda=U,$ where $U=(-\infty,a]$ and $L=\emptyset$;
\item $\Lambda=L,$ where $L=[b,\infty)$ and $U=\emptyset$;
\item $\Lambda=(a,b)$ where $L=(-\infty,a]$ and $U=[b,\infty)$.
\end{enumerate}
If $(\xi,\eta)$ has $L\cap U=\emptyset$ and does not have $U$ and
$L$ satisfying one of (1)-(4), then no absorbing sets exist.
\end{cor}

We examine two striking cases of $L$ and $U$ structure, and state
exact conditions on the characteristic triplet of $(\xi,\eta)$ for
such behaviour. Note that similar conditions can be found for each
of the other $L$ and $U$ structures stated in Theorem
\ref{combinations proposition}, however, the statements are longer
and unwieldy.
\begin{prop}\label{convergent divergent proposition}
Suppose $L\cap U=\emptyset.$ Then $U=(-\infty,a]$ and $L=[b,\infty)$
for $-\infty<a<b<\infty$ iff $(\xi,\eta)$ is of finite variation and
the following hold:
\begin{itemize}
\item There is no Brownian component $(\Sigma_{\xi,\eta}=0);$
\item The drift of $\xi$ is non-negative $(d_\xi\ge 0);$
\item The \Levy measure satisfies $\Pi_{\xi,\eta}(A_3)=\Pi_{\xi,\eta}(A_4)=0,$ $\theta_1'>-\infty,$ and $\theta_2<\infty.$
\end{itemize}
If these conditions hold then $\xi$ is a subordinator and, for any
$V_0=z\in\mathbb{R},$ $\lim_{t\rightarrow\infty}|V_t|=\infty$ a.s.

Similarly $L=(-\infty,a]$ and $U=[b,\infty)$ for
$-\infty<a<b<\infty$ iff $(\xi,\eta)$ is of finite variation and the
following hold:
\begin{itemize}
\item There is no Brownian component $(\Sigma_{\xi,\eta}=0);$
\item The drift of $\xi$ is non-positive $(d_\xi\le 0);$
\item The \Levy measure satisfies $\Pi_{\xi,\eta}(A_1)=\Pi_{\xi,\eta}(A_2)=0,$
$\theta_4'<\infty$ and $\theta_3>-\infty.$
\end{itemize}
If these conditions hold then $-\xi$ is a subordinator, and $V$ is
strictly stationary and converges in distribution as
$t\rightarrow\infty$ to a random variable $V_\infty$ supported on
$(a,b).$
\end{prop}

We now present a theorem describing the relationship between the
sets $L$ and $U,$ and the upper and lower bounds of the limit random
variable $Z_\infty$ of $Z_t$ as $t\rightarrow\infty.$

\begin{thm}\label{linking theorem} Let $a,b\in\mathbb{R}$ and suppose $Z_t\rightarrow Z_\infty$  a.s. as
$t\rightarrow\infty,$ where $Z_\infty$ is a finite random variable.
If, for all $c\in\mathbb{R},$ the degenerate case (\ref{degenerate
case equation}) does not hold, then $a\le\sup U$ iff $Z_\infty<-a$
a.s., whilst $b\ge\inf L$ iff $Z_\infty>-b$ a.s. Further, $-\sup
U=\inf\{u\in\mathbb{R}|Z_\infty<u~\mathrm{a.s.}\}$ and $-\inf
L=\sup\{u\in\mathbb{R}|Z_\infty>u~\mathrm{a.s.}\}.$ Alternatively,
if there exists $c\in\mathbb{R}$ such that equation (\ref{degenerate
case equation}) holds, then $Z_\infty=-c$ a.s. and $\inf L=\sup
U=c.$
\end{thm}

The next theorem presents results on certain ruin which occur when
$L$ and $U$ are of a particular structure.

\begin{thm}\label{minor certain ruin theorem} Suppose that $L\cap
U=\emptyset.$ Then the following statements hold:
\begin{enumerate}
\item If $\sup U\ge 0$ and $L\cap[0,\sup U]=\emptyset,$ then
$\psi(z)=1$ for all $z\le\sup U;$
\item If $\sup L\ge 0$ and $U\cap[0,\sup L]=\emptyset,$ then
$0<\psi(z)<1$ for all $0\le z<\inf L.$ If $\sup L\ge 0$ and
$U\cap[0,\sup L]\neq\emptyset,$ then $\psi(z)<1$ for all $z>\sup U.$
\end{enumerate}
\end{thm}

Note that in statement (2) above, when $\sup L\ge0$ and $L\cap
U\neq\emptyset,$ Theorem \ref{combinations proposition} ensures that
$\sup U<\inf L,$ and statement (1) above ensures that $\psi(z)=1$
for all $z\le \sup U.$ Also, by definition of $L,$ $\psi(z)=0$
whenever $z\ge\inf L.$

We now present a major theorem which utilises Theorems
\ref{combinations proposition}, \ref{linking theorem} and \ref{minor
certain ruin theorem}, and is the major tool in proving Theorems
\ref{paulsen theorem 1} and \ref{paulsen theorem 2}. For the
non-degenerate case, and for $(\xi,\eta)$ which satisfies various
asymptotic and stability criteria, this theorem presents iff
conditions for certain ruin, stated in terms of $L$ and $U$
structure. In particular, it completely describes the $L$ and $U$
structures for which certain ruin occurs.

\begin{thm}\label{ruin prob theorem for L and U} Suppose $L\cap U=\emptyset.$
\begin{enumerate}
\item Suppose $\lim_{t\rightarrow\infty}\xi_t=-\infty$ a.s. and
$I_{-\xi,K^{\xi,\eta}}<\infty.$ There exists $z\ge 0$ such that
$\psi(z)<1$
 iff $L\cap[0,\infty)\neq \emptyset.$ If this occurs then
$0<\psi(z)<1$ for all $0\le z<\inf L,$ $\psi(z)=0$ for all $z\ge\inf
L,$ and one of the following must hold:
\begin{enumerate}
\item $L=[a,b]$ and $U=\emptyset,$ where $-\infty\le a\le b<\infty,$ and $b\ge 0;$
\item $L=(-\infty,a]$ and $U=[b,\infty)$ where $0\le a<b<\infty.$
\end{enumerate}
\item Suppose $E(\xi_1)=0,$ $E(e^{|\xi_1|})<\infty$ and there exist $p,q>1$ with
$1/p+1/q=1$ such that $E\left(e^{-p\xi_1}\right)<\infty$ and
$E\left(|\eta_1|^q\right)<\infty.$ There exists $z\ge 0$ such that
$\psi(z)<1$ iff $L\cap[0,\infty)\neq \emptyset.$ If this occurs then
$L=[a,b]$ and $U=\emptyset,$ where $-\infty< a\le b<\infty$ and
$b\ge 0,$ in which case $0<\psi(z)<1$ for all $0\le z<a$ and
$\psi(z)=0$ for all $z\ge a;$

\item Suppose
$\lim_{t\rightarrow\infty}\xi_t=\infty$ a.s. and
$I_{\xi,\eta}<\infty.$ There exists $z\ge 0$ such that $\psi(z)=1$
iff $U\cap[0,\infty)\neq\emptyset.$ If this occurs then one of the
following must hold:
\begin{enumerate}
\item [(c)] $U=[a,b]$ and $L=\emptyset,$ where $-\infty\le a\le b<\infty$ and $b\ge
0,$ in which case $\psi(z)=1$ for all $z\le b$ and $0<\psi(z)<1$ for
all $z>b;$

\item [(d)] $U=(-\infty,a]$ and $L=[b,\infty)$ where $0\le a<
b<\infty,$ in which case $\psi(z)=1$ for all $z\le a,$ $0<\psi(z)<1$
for all $a<z<b$ and $\psi(z)=0$ for all $z\ge b.$

\end{enumerate}
\end{enumerate}
\end{thm}

\begin{rem} \rm{The characteristic triplet conditions which equate to the iff  result in statement (3)
above, are given in Theorem \ref{paulsen theorem 2}, and are
obtained using the forthcoming Proposition \ref{upper bound
theorem}. Further, exact characteristic triplet conditions for the
structure $U=(-\infty,a]$ and $L=[b,\infty)$ in case (d) above, are
given in Proposition \ref{convergent divergent proposition}. }
\end{rem}

\subsection{Technical results on the upper and lower bounds}

We present a series of important technical propositions on $\delta,$
$L,$ $\Upsilon$ and $U.$ As well as being of independent interest,
they are essential in proving the previously stated theorems. The
first proposition is obtained by combining and restating parts of
Proposition 6, Theorem 7 and Theorem 9 of \cite{BankovskySly08}, and
no proof is given. When put into this form the proposition
completely describes the relationship between the \Levy measure of
$(\xi,\eta)$ and the lower bound function $\delta.$ We recall some
notation from \cite{BankovskySly08}. For $A_i$ as in Section \ref{ruin probability section}, define
$A_i^u:=\left\{(x,y)\in A_i:y-u(e^{-x}-1)<0\right\}.$ For $u\le 0$
define
\[
\theta_1:= \left\{ \begin{array}{ll}
\sup\left\{u\le 0:\Pi_{\xi,\eta}(A_1^u)>0\right\} & \\
-\infty~~~~\textrm{if~}\Pi_{\xi,\eta}(A_1\setminus A_4)=0, &
\end{array} \right.\,
\theta_3:= \left\{ \begin{array}{ll}
\inf\left\{u\le 0:\Pi_{\xi,\eta}(A_3^u)>0\right\} &\\
0~~~~\textrm{if~}\Pi_{\xi,\eta}(A_3\setminus A_4)=0, &
\end{array} \right.\,
\]
and for $u\ge 0$ define
\[
\theta_2:= \left\{ \begin{array}{ll}
\sup\left\{u\ge 0:\Pi_{\xi,\eta}(A_2^u)>0\right\} & \\
0~~~~\textrm{if~}\Pi_{\xi,\eta}(A_2\setminus A_1)=0, &
\end{array} \right.\,
\theta_4:= \left\{ \begin{array}{ll}
\inf\left\{u\ge 0:\Pi_{\xi,\eta}(A_4^u)>0\right\} &\\
\infty~~~~\textrm{if~}\Pi_{\xi,\eta}(A_4\setminus A_1)=0. &
\end{array} \right.\,
\]
Throughout, let $W$ be the \Levy process such that
$e^{-\xi_t}=\epsilon(W)_t.$

\begin{prop}[lower bound]\label{lower bound theorem} The following statements are equivalent:
\begin{enumerate}
\item The lower bound $\delta(z)>-\infty$ for some $z\in\mathbb{R}$;
\item There exists $u\in\mathbb{R}$ such that $\delta(u)=u$;
\item There exists $u\in\mathbb{R}$ such that the \Levy process $\eta-uW$ is a subordinator.
\end{enumerate}
Statements (2) and (3) hold for a particular value $u\neq0$ iff the
following three conditions are satisfied: (i) the Gaussian
covariance matrix satisfies equation (\ref{covariance matrix
equation}); (ii) one of the following is true:
\begin{enumerate}
\item[(a)] $\Pi_{\xi,\eta}(A_3)=0,$ $\Pi_{\xi,\eta}(A_2)\neq0,$ $\theta_2\le\theta_4$ and $u\in[\theta_2,\theta_4];$
\item[(b)] $\Pi_{\xi,\eta}(A_2)=0,$ $\Pi_{\xi,\eta}(A_3)\neq0,$ $\theta_1\le\theta_3$ and $u\in[\theta_1,\theta_3];$
\item[(c)] $\Pi_{\xi,\eta}(A_3)=\Pi_{\xi,\eta}(A_2)=0$ and $u\in[\theta_1,\theta_4];$
\end{enumerate}
 and, (iii), in addition, $u$ satisfies $g(u)\ge 0$ for the
 function g in equation (\ref{finite drift equation}).

\end{prop}
From the definition of $L$ it is an immediate corollary, firstly,
that $L=\emptyset$ iff none of conditions (1)-(3) of Proposition
\ref{lower bound theorem} hold, and secondly, that $\eta$ is a
subordinator iff $0\in L.$ The next proposition adds further
information concerning $L.$ Most importantly, it shows that the set
$L$ is always connected, and gives concrete values for the
endpoints.

\begin{prop}\label{endpoints of L proposition}
If $\sigma^2_\xi\neq 0$ and any of conditions (1)-(3) of Proposition
\ref{lower bound theorem} hold, then
$L=\{-\frac{\sigma_{\xi,\eta}}{\sigma^2_\xi}\}.$ If $\sigma^2_\xi=
0$  and any of (1)-(3) hold, then $\sigma^2_\eta=0$ and one of the
following holds:
\begin{itemize}
\item $\eta$ is a subordinator and condition (ii) of Proposition
\ref{lower bound theorem} does not hold for any $u\neq 0$, in which
case $L=\{0\};$
\item Condition (ii) is satisfied for some $u\neq0$, in which case there exists $-\infty\le a\le b\le \infty$ such that $L=[a,b].$
\end{itemize}
In the latter case, if condition (a) of Proposition \ref{lower bound
theorem} holds then $0\le a=\max\{\theta_2,m_1\}$ and
$b=\min\{\theta_4,m_2\}$ for $m_1:=\inf\{u\in\mathbb{R}:g(u)\ge 0\}$
and $m_2:=\sup\{u\in\mathbb{R}:g(u)\ge 0\}.$ If (b) holds then
$a=\max\{\theta_1,m_1\}$ and $b=\min\{\theta_3,m_2\}\le 0.$ If (c)
holds then $a=\max\{\theta_1,m_1\}$ and $b=\min\{\theta_4,m_2\}.$
\end{prop}

Define $L^*$ to be the set of starting values on which the GOU has
no negative jumps, namely
\[L^*:=\left\{u\in\mathbb{R}:\forall t>0~P\left(\Delta
V_t<0|V_{t-}=u\right)=0\right\}.
\]
It is a consequence of Proposition \ref{description of L
proposition} that $L\subseteq L^*.$ The next proposition describes
$L^*.$  In particular, it shows that the set $L^*$ is always
connected, and gives concrete values for the endpoints. It also
shows that whenever $V_{t-}>\sup L^*$ and a negative jump $\Delta
V_t$ occurs, then the jump cannot be so negative as to cause $V_t\le
\sup L^*.$ Thus, $L^*$ acts as a barrier for negative jumps of $V.$

\begin{prop}\label{endpoints of L star proposition}
\begin{enumerate}

\item If $L^*\neq\emptyset$ then, for any $t\ge 0,$
$V_{t-}>\sup L^*$ implies $V_t>\sup L^*$ a.s.;

\item $L^*=\{u\in\mathbb{R}:\eta-uW\mathrm{~has~no~negative~jumps}\};$

\item  $L^*\neq \emptyset$ iff condition (ii) of Proposition
\ref{lower bound theorem} is satisfied for some $u\neq 0,$ or $\eta$
has no negative jumps;

\item $L^*=\{0\}$ iff $\eta$ has no negative jumps and condition (ii) does not hold for
any $u\neq 0;$

\item If condition (ii) of Proposition \ref{lower bound theorem} holds
for some $u\neq 0$ then $L^*=[\theta_2,\theta_4],$
$[\theta_1,\theta_3]$ or $[\theta_1,\theta_4],$ corresponding to
conditions (a), (b) or (c) of Proposition \ref{lower bound theorem}.

\end{enumerate}\end{prop}

\begin{rem}\label{lower bound remarks}
\begin{enumerate}
\rm{
\item If $(\xi,\eta)$ is an infinite variation \Levy process then,
as noted in Section \ref{intro section},
$\int_{\{x^2+y^2<1\}}|(x,y)|\Pi_{\xi,\eta}(\ud (x,y))=\infty.$ Thus,
it may be the case that for a particular $u\in\mathbb{R}$ the
integral $\int_{\{x^2+y^2<1\}}(ux+y)\Pi_{\xi,\eta}(\ud (x,y)),$ and
hence the function $g(u)$ in (\ref{finite drift equation}), may not
exist as a well-defined member of the extended real numbers.
However, it is a consequence of the proof of Theorem 9 in
\cite{BankovskySly08}, that if $u\in L^*$ then $g(u)$ is a well
defined member of the extended reals, and $g(u)\in[-\infty,\infty).$
Under such conditions, it is also shown that
\[\Pi_{\xi,\eta}\left(\{y-u(e^{-x}-1)< 0\}\right)=0
\]
and so the domain of integration for the integral component of $g$
can be decreased to $\{x^2+y^2<1\}\cap\{y-u(e^{-x}-1)\ge 0\}$.

\item Note that $g$ is a linear function on $\mathbb{R}$ iff
the \Levy measure of $(\xi,\eta)$ is of finite variation, namely
\[\int_{\{x^2+y^2<1\}}|(x,y)|\Pi_{\xi,\eta}(\ud(x,y))<\infty.\]
 In this case the drift vector $(d_\xi,d_\eta)$ is finite, and we
 can write
\begin{eqnarray}\label{linear g
equation}g(u)&=&\gamma_\eta-\int_{(-1,1)}y\Pi_\eta(\ud
y)+u\left(\gamma_\xi-\frac{1}{2}\sigma_\xi^2-\int_{(-1,1)}x\Pi_\xi(\ud
x)\right)\nonumber\\
&=&d_\eta+u\left(d_\xi-\frac{1}{2}\sigma_\xi^2\right),
\end{eqnarray}
where the first equality follows by converting
$(\tilde{\gamma}_\xi,\tilde{\gamma}_\eta)$ to
$(\gamma_\xi,\gamma_\eta)$ using equation (\ref{first 2 dim to 1 dim
equation}) and the symmetric version for $\eta$, and the second
equality follows by converting $(\gamma_\xi,\gamma_\eta)$ to
$(d_\xi,d_\eta)$ using equation (\ref{drift equation}). It will be a
consequence of the proof of Proposition \ref{endpoints of L
proposition}, that if $a,b\in L$ and $a\neq b$ then $g$ is a linear
function on $\mathbb{R}.$

\item In Section \ref{intro section} we stated exact conditions for a \Levy process to be a subordinator.
When $u\neq 0$ the \Levy measure conditions in Proposition
\ref{lower bound theorem} are exactly the requirements for $\eta-uW$
to be a subordinator. Equation (\ref{covariance matrix equation}) is
equivalent to the condition $\sigma_{\eta-uW}=0.$ The requirement
that one of the conditions (a), (b) and (c) holds is equivalent to
the requirement that there exists $u\neq 0$ such that
$\Pi_{\eta-uW}((-\infty,0))=0.$ Note that this implies that
$L^*\setminus\{0\}$ is precisely the set of all $u\neq 0$ such
$\eta-uW$ has no negative jumps. Finally, if $u\in L^*$ then
$g(u)=d_{\eta-uW},$ and hence condition (\ref{finite drift
equation}) is equivalent to the requirement that $\eta-uW$ has
positive drift. The fact that $\eta-uW$ is of finite variation
actually follows from the two conditions
$\Pi_{\eta-uW}((-\infty,0))=0$ and $d_{\eta-uW}\ge 0.$ To see this,
note that when $\Pi_{\eta-uW}((-\infty,0))=0,$ the equation
(\ref{drift equation}) simplifies to
\[d_{\eta-uW}=\gamma_{\eta-uW}-\int_{(0,1)}x\Pi_{\eta-uW}(\ud x)
\] and hence $d_{\eta-uW}$ is a member of the extended
reals regardless of whether $\eta-uW$ is finite variation. In
particular, $d_{\eta-uW}\in[-\infty,\infty),$ and
$d_{\eta-uW}=-\infty$ iff $\int_{(0,1)}x\Pi_{\eta-uW}(\ud x)=\infty$
which occurs iff $\eta-uW$ is infinite variation. }\end{enumerate}
\end{rem}

Although the situation is symmetric, we explicitly state the
parallel version for $U$ and $\Upsilon,$ to Proposition \ref{lower
bound theorem}. No proof is given. We state the parallel result
explicitly because some of the statements are not obvious, and we
need to use them for Theorem \ref{paulsen theorem 2}. Also, we will
need to combine them with the statements for $L$ and $\delta$ in
order to prove Theorem \ref{combinations proposition}, \ref{minor
certain ruin theorem} and \ref{ruin prob theorem for L and U}.  If
we define
\[U^*:=\left\{u\in\mathbb{R}:\forall t>0~P\left(\Delta
V_t>0|V_{t-}=u\right)=0\right\},
\] then the symmetric versions of Proposition \ref{endpoints of L
proposition}, Proposition \ref{endpoints of L star proposition} and
Remark \ref{lower bound remarks} also hold. We will need to use
these results, however the parallels are obvious in this case, so we
do not state them explicitly.

\begin{prop}[upper bound]\label{upper bound theorem} The following are equivalent:
\begin{enumerate}
\item The upper bound $\Upsilon(z)<\infty$ for some $z\in\mathbb{R}$;
\item There exists $u\in\mathbb{R}$ such that $\Upsilon(u)=u$;
\item There exists $u\in\mathbb{R}$ such that the \Levy process $-(\eta-uW)$ is a
subordinator.
\end{enumerate}
Statements (2) and (3) hold for a particular value $u\neq 0$ iff the
following three conditions are satisfied: (i) the Gaussian
covariance matrix satisfies equation (\ref{covariance matrix
equation}); (ii) one of the following is true:
\begin{enumerate}
\item[(a)] $\Pi_{\xi,\eta}(A_1)=0,$ $\Pi_{\xi,\eta}(A_4)\neq0,$ $\theta_4'\le\theta_2'$ and $u\in[\theta_4',\theta_2'];$
\item [(b)]$\Pi_{\xi,\eta}(A_4)=0,$ $\Pi_{\xi,\eta}(A_1)\neq0,$ $\theta_3'\le\theta_1'$ and $u\in[\theta_3',\theta_1'];$
\item[(c)] $\Pi_{\xi,\eta}(A_1)=\Pi_{\xi,\eta}(A_4)=0$ and $u\in[\theta_3',\theta_2'];$
\end{enumerate}
and,(iii), in addition, $u$ satisfies $g(u)\le0$ for the function
$g$ in equation (\ref{finite drift equation}).
\end{prop}

\begin{rem}\label{upper bound remarks} \rm{Symmetric statements to those for $L$ and $L^*$ in Remark \ref{lower bound
remarks}, hold for $U$ and $U^*$. The following remarks relate to
the combination of $L$ and $U,$ and $L^*$ and $U^*.$\rm}
\begin{enumerate}
\rm{
\item Parallel to 1 and 2 of Remark \ref{lower bound remarks},
whenever $u\in U^*,$ $g(u)$ from (\ref{finite drift equation}) is a
well-defined member of the extended reals,
$g(u)\in(-\infty,\infty],$ and $-g(u)=d_{-(\eta-uW)}.$ Since
$d_{-(\eta-uW)}=-d_{\eta-uW},$ we know that if $u\in U^*\cup L^*$
then $g(u)$ is a well-defined member of the extended reals and
$g(u)=d_{\eta-uW}.$

\item If $a\in L,$ $b\in U$ and $a\neq b$ then $g$ is linear and
$(\xi,\eta)$ is finite variation. This statement is proved easily
using similar arguments to those in the proof of Proposition
\ref{endpoints of L proposition}.

\rm}
\end{enumerate}
\end{rem}

We state a proposition, describing the possible combinations of
$L^*$ and $U^*,$ which will be essential for proving Theorem
\ref{combinations proposition}.

\begin{prop}\label{combinations of L star and U star lemma}
The following statements hold for $L^*,$ and the symmetric
statements hold for $U^*:$
\begin{enumerate}
\item $L^*=\mathbb{R}$ then $U^*=\emptyset$ or $U^*=\mathbb{R};$
\item If $L^*=[a,b]$ for some $-\infty<a\le b<\infty,$ then
$U^*=\emptyset$ or $U^*=L^*=\{a\}=\{b\};$
\item If $L^*=[b,\infty)$ for some $b\in\mathbb{R},$ then $U^*=\emptyset$ or
$U^*=(-\infty,a]$ for some $-\infty<a\le b<\infty;$
\item If
$L^*=(-\infty,a]$ for some $a\in\mathbb{R},$ then $U^*=\emptyset$ or
$U^*=[b,\infty)$ for some $-\infty<a\le b<\infty.$
\end{enumerate}
\end{prop}

We end the section with two lemmas. No proof will be given. The
first follows by considering the definitions of $\theta_i$ and
$\theta_i'.$ It will be used several times as a calculation tool.
The second gives conditions on the \Levy measure of $\xi$ and $\eta$
which ensure that $\sup_{0\le t\le 1}|Z_t|$ has finite mean. It will
be needed to prove statement (2) of Theorem \ref{paulsen theorem 1}.
The proof is similar to that of Lemma 11 in \cite{BankovskySly08}
and uses the Burkholder-Davis-Gundy inequalities, and various Doob's
inequalities.

\begin{lem}\label{inequalities proposition}
\begin{enumerate}
\item If $\Pi_{\xi,\eta}(A_1)\neq0$ then $\theta_1'\le\theta_1\le0;$
\item If $\Pi_{\xi,\eta}(A_2)\neq0$ then $0\le\theta_2'\le\theta_2;$
\item If $\Pi_{\xi,\eta}(A_3)\neq0$ then $\theta_3\le\theta_3'\le0;$
\item If $\Pi_{\xi,\eta}(A_4)\neq0$ then $0\le\theta_4\le\theta_4'.$
\end{enumerate}
Further:
\begin{enumerate}
\item[(a)] $\Pi_{\xi,\eta}(A_1)=0$ iff $\theta_1=-\infty$ and $\theta_1'=0;$
\item[(b)] $\Pi_{\xi,\eta}(A_2)=0$ iff $\theta_2=0$ and $\theta_2'=\infty;$
\item[(c)] $\Pi_{\xi,\eta}(A_3)=0$ iff $\theta_3=0$ and $\theta_3'=-\infty;$
\item[(d)] $\Pi_{\xi,\eta}(A_4)=0$ iff $\theta_4=\infty$ and $\theta_4'=0.$
\end{enumerate}
\end{lem}

\begin{lem}\label{sup proposition} Suppose there exist
$r>0$ and $p,q>1$ with $1/p+1/q=1$ such that $E\left(e^{-\max \{1,r
\}p\xi_1}\right)<\infty$ and $E\left(|\eta_1| ^{\max\{1,r
\}q}\right)<\infty.$ Then
\begin{equation}E\left(\sup_{0\le t\le 1}\left|\int_0^t e^{-\xi_{s-}}\ud\eta_s\right|^{\max
\{1,r \}}\right)<\infty .
\end{equation}
\end{lem}

\section{Proofs and Examples}

The proofs are presented in mathematically chronological order
rather than the order in which the statements of the results are
presented. For all proofs, except the proof of Proposition
\ref{degenerate proposition}, we assume that neither $\xi$ nor
$\eta$ are zero.

\begin{pf}[Proposition \ref{endpoints of L star proposition}]
We prove statements (1), (2) and (3). The proof of statements (4)
and (5) follows trivially from the proof of statements (2) and (3).

(1) Suppose $L^*=\emptyset.$ Assume that condition (a) of
Proposition \ref{lower bound theorem} holds and
$L^*=[\theta_2,\theta_4].$ If condition (b) or (c) of Proposition
\ref{lower bound theorem} holds then the proof is similar. We use
the following reformulation of equation (\ref{jump equation 1}):
\begin{equation}\label{jump equation 2}
\Delta V_t=(e^{\Delta\xi_t}-1)V_{t-}+e^{\Delta\xi_t}\Delta\eta_t.
\end{equation}
Suppose $V_{t-}>\theta_4.$ It follows immediately from the
definitions of $\theta_4$ and $A_4^u,$ and from equation (\ref{jump
equation 2}), that there exists $(x,y)\in A_4^{V_{t-}}$ such that
$(e^x-1)\theta_4+e^xy\ge 0$ and $(e^x-1)V_{t-}+e^xy<0.$ Thus,
\begin{eqnarray*}V_t&=&V_{t-}+(e^x-1)V_{t-}+e^xy\\
&=&V_{t-}+(e^x-1)(V_{t-}-\theta_4)+(e^x-1)\theta_4+e^xy\\
&\ge&V_{t-}+(e^x-1)(V_{t-}-\theta_4)>\theta_4.
\end{eqnarray*}
(2) It is a consequence of Proposition 6 in
\cite{BankovskySly08} that
\[\Delta(\eta_t-uW_t)=\Delta\eta_t-u\left(e^{-\Delta\xi_t}-1\right).\]
Thus,  equation (\ref{jump equation 1}) implies that whenever
$V_{t-}=u,$ a jump $(\Delta \xi_t,\Delta\eta_t)$ causes a negative
jump $\Delta V_t$ iff $\Delta(\eta_t-uW_t)$ is negative. Hence $L^*$
is precisely the set of all $u$ such that $\eta_t-uW_t$ has no
negative jumps.

(3) By (1) above, $L^*\neq\emptyset$ iff $\eta-uW$ has no negative
jumps. If $u=0,$ this occurs iff $\eta$ has no negative jumps. If
$u\neq0,$ it is noted in point (3) of Remark \ref{lower bound
remarks}, that this occurs iff $u\neq0$ satisfies condition (ii) of
Proposition \ref{lower bound theorem}. \halmos
\end{pf}

\begin{pf}[Proposition \ref{endpoints of L proposition}]
Assume that $\sigma^2_\xi\neq 0$ and statements (1)-(3) of
Proposition \ref{lower bound theorem} hold for some $u\neq 0.$ Then
equation (\ref{covariance matrix equation}) must hold for $u,$ which
implies that $u=-\frac{\sigma_{\xi,\eta}}{\sigma^2_\xi},$ and hence
is the unique non-zero number satisfying statements (1)-(3) of
Proposition \ref{lower bound theorem}. Since
$-\frac{\sigma_{\xi,\eta}}{\sigma^2_\xi}$ satisfies condition (2),
$L=\{-\frac{\sigma_{\xi,\eta}}{\sigma^2_\xi}\}$ by definition.

Now assume that $\sigma^2_\xi\neq 0$ and statements (1)-(3)
of Proposition \ref{lower bound theorem} hold for $u=0.$ By
statement (2), $0\in L.$ By statement (3), $\eta$ is a subordinator,
and hence $\sigma^2_\eta=\sigma_{\xi,\eta}=0.$ Thus, by the above,
no non-zero number can satisfy statements (1)-(3), and so
$L=\{0\}=\{-\frac{\sigma_{\xi,\eta}}{\sigma^2_\xi}\}.$

Now assume that $\sigma^2_\xi=0.$ If statements (1)-(3) of
Proposition \ref{lower bound theorem} hold for $u=0$ then $\eta$ is
a subordinator by statement (3) and hence $\sigma^2_\eta=0.$
Alternatively, If statements (1)-(3) of Proposition \ref{lower bound
theorem} hold for some $u\neq0$ then equation (\ref{covariance
matrix equation}) must hold for $u,$ which implies that
$\sigma^2_\eta=u^2\sigma^2_\xi,$ and so $\sigma^2_\eta=0.$

Now assume that $\sigma^2_\xi=0$ and condition (ii) of Proposition
\ref{lower bound theorem} does not hold for any $u\neq 0.$ This
immediately implies that
$L\cap(\mathbb{R}\setminus\{0\})=\emptyset.$ If, further, $\eta$ is
a subordinator, then $0\in L,$ and hence $L=\{0\}.$

Now assume that $\sigma^2_\xi=0$ and condition (ii) of Proposition
\ref{lower bound theorem} holds for some $u\neq 0.$ This occurs
precisely when one of conditions (a), (b) or (c) of Proposition
\ref{lower bound theorem} holds, and equation (\ref{finite drift
equation}) holds. Thus, $\inf L=a$ and $\sup L=b$ for the values of
$a$ and $b$ given in the proposition statement. It remains to prove
that the set $L$ is connected. Since $L^*$ is connected, this occurs
iff $\{u\in\mathbb{R}:g(u)\ge 0\}$ is connected, which follows from
the analysis below.

As noted in point (1) of Remark \ref{lower bound remarks}, whenever
$u\in L^*$ we know $g(u)\in[-\infty,\infty).$ There are
three possibilities for behaviour of $g$ on $L^*.$ Firstly, it may
be that $g(u)=-\infty$ for all $u\in L^*.$ Secondly there may exist
$v\in L^*$ such that $g(v)$ is finite and $g(u)=-\infty$ for all
$u\in L^*$ with $u\neq v.$ We show that the only other possibility
is that $g$ is linear on $\mathbb{R}.$ Suppose there exists
$u_1,u_2\in L^*$ with $u_1\neq u_2,$ such that $g(u_1)$ and $g(u_2)$
are both finite. Then
\[g(u_1)-g(u_2)=\left(\tilde{\gamma_\xi}-\frac{1}{2}\sigma_\xi^2-\int_{\{x^2+y^2<1\}}x\Pi_{\xi,\eta}(\ud(x,y))\right)(u_1-u_2)
\]
is finite, which implies that
$\int_{\{x^2+y^2<1\}}x\Pi_{\xi,\eta}(\ud(x,y))$ exists, and is
finite. Since $g(u_1)$ is finite, this implies that
$\int_{\{x^2+y^2<1\}}y\Pi_{\xi,\eta}(\ud(x,y))$ exists and is
finite. Thus, $g$ is a linear function on $\mathbb{R}.$ \halmos
\end{pf}

\begin{pf}[Proposition \ref{description of L proposition}]
It is an immediate consequence  of Proposition \ref{lower bound
theorem} that $\delta(\delta(z))=\delta(z)$ and
\begin{equation}\label{linking delta and L equation}\delta(z)=\sup\{u\le z:\delta(u)=u\}.
\end{equation}
Now the first statement of Proposition \ref{description of L
proposition} follows immediately from (\ref{linking delta and L
equation}). To prove the second statement, assume $z<\inf L.$
Suppose $-\infty<m:=\delta(z).$ Since $\delta(z)\le z,$ we have
$-\infty<m\le z<\inf L.$ However, equation (\ref{linking delta and L
equation}) implies that $m\in L,$ which gives a contradiction. Hence
$\delta(z)=-\infty.$ The third and fourth statements follow
immediately from the definitions of $\delta$ and $L.$ \halmos
\end{pf}

\begin{pf}[Proposition \ref{degenerate proposition}] Assume $L=\mathbb{R}.$ This implies, using Proposition \ref{lower bound theorem}
and point (2) of Remark \ref{lower bound remarks}, that
$\Sigma_{\xi,\eta}=0$ and $g$ is linear. Further, it must be the
case that $\Pi_{\xi,\eta}(A_3)=\Pi_{\xi,\eta}(A_2)=0$ and
$L^*=[\theta_1,\theta_4]=(-\infty,\infty).$ Now $\theta_1=-\infty$
iff $\Pi_{\xi,\eta}\left((0,\infty)\times[0,\infty)\right),$ whilst
$\theta_4=-\infty$ iff
$\Pi_{\xi,\eta}\left((-\infty,0)\times[0,\infty)\right)=0.$ Hence
$\xi$ can have no jumps and $\eta$ can only have positive jumps. By
Proposition \ref{lower bound theorem}, $g(u)\ge0$ on $\mathbb{R}.$
Since $g(u)=d_\eta+ud_\xi,$ this implies that $d_\xi=0$ and
$d_\eta\ge0,$ thus proving one direction of the first claim. The
converse is trivial since $V$ simplifies to $V_t=z+\eta_t.$ The
proof of the second claim is similar. The third claim follows
immediately from the first two. \halmos
\end{pf}

\begin{pf}[Proposition \ref{combinations of L star and U star lemma}]
We prove statements (1), (2) and (3). The proof of statement (4) is
symmetrical to the proof of statement (3).

(1) Assume $L^*=\mathbb{R}.$ Then condition (c) of Proposition
\ref{lower bound theorem} must hold, and so  $\Pi_{\xi,\eta}(A_2)=
\Pi_{\xi,\eta}(A_3)=0,$ and $L^*=[\theta_1,\theta_4].$ Since
$\theta_1=-\infty$ and $\theta_4=\infty,$ it must be that
$\Pi_{\xi,\eta}(A_1\setminus A_4)=0$ and
$\Pi_{\xi,\eta}(A_4\setminus A_1)=0,$ respectively. Thus, if
$\Pi_{\xi,\eta}(A_1\cap A_4)=0$ then
$\Pi_{\xi,\eta}(\mathbb{R}^2)=0,$ in which case condition (c) of
Proposition \ref{upper bound theorem} holds, and $U^*=\mathbb{R}$.
Alternatively, if $\Pi_{\xi,\eta}(A_1\cap A_4)\neq0$ then $\eta$ has
positive jumps and so $0\not\in U^*,$ and (ii) of Proposition
\ref{upper bound theorem} cannot hold. Hence $U^*=\emptyset.$

(2) Assume $L^*=[a,b]$ for some $-\infty<a\le b<\infty.$ There are
four ways in which this is possible, namely, when (a), (b) or (c) of
Proposition \ref{lower bound theorem} hold, or when $L^*=\{0\}.$ For
each case we show $U^*=\emptyset$ or $U^*=L^*=\{a\}=\{b\}.$

Suppose first that condition (a) of Proposition \ref{lower bound
theorem} holds, and $U^*\neq\emptyset.$ The case in which condition
(b) holds and $U^*\neq\emptyset,$ is symmetric. Propositions
\ref{lower bound theorem} and \ref{endpoints of L star proposition}
imply that $\Pi_{\xi,\eta}(A_3)=0,$ $\Pi_{\xi,\eta}(A_2)\neq 0,$
$\theta_2\le \theta_4$ and $L^*=[\theta_2,\theta_4].$ Since
$\theta_4<\infty,$ it must be that $\Pi_{\xi,\eta}(A_4\setminus
A_1)\neq 0.$ Since $\Pi_{\xi,\eta}(A_3)=0,$ this implies that
$-\eta$ is not a subordinator, and so $0\not\in U^*.$ Thus, since we
have assumed that $U^*\neq\emptyset,$ it must be that condition (a)
of Proposition \ref{upper bound theorem} holds, and so
$\Pi_{\xi,\eta}(A_1)=0,$ $\theta_4'\le\theta_2',$ and
$U^*=[\theta_4',\theta_2'].$ However, statements (2) and (4) of
Lemma \ref{inequalities proposition} state that
$\theta_2'\le\theta_2$ and $\theta_4\le\theta_4'.$ Hence
$\theta_2'=\theta_2=\theta_4=\theta_4'.$

Now suppose that condition (c) of Proposition \ref{lower bound
theorem} holds. Then $\Pi_{\xi,\eta}(A_2)=\Pi_{\xi,\eta}(A_3)=0,$
and $L^*=[\theta_1,\theta_4].$ Since $\theta_4<\infty$ and
$\theta_1>-\infty$ it must be that $\Pi_{\xi,\eta}(A_4\setminus
A_1)\neq 0$ and $\Pi_{\xi,\eta}(A_1\setminus A_4)\neq 0,$
respectively. Hence condition (ii) of Proposition \ref{upper bound
theorem} cannot hold, and so $U^*\setminus\{0\}=\emptyset.$ Further,
$-\eta$ is not a subordinator, and so $U^*=\emptyset.$

Now suppose $L^*=\{0\},$ and $U^*\neq\emptyset.$  By statement (4)
of Proposition \ref{endpoints of L star proposition}, $L^*=\{0\}$
iff $\eta$ has no negative jumps and at the same time
$\Pi_{\xi,\eta}(A_3\cap A_4)\neq 0$ and $\Pi_{\xi,\eta}(A_2\cap
A_1)\neq 0.$ Hence, condition (ii) of Proposition \ref{upper bound
theorem} fails to hold, which implies $U^*\setminus\{0\}=\emptyset.$
Thus, since $U^*\neq\emptyset,$ it must be that $U^*=L^*=\{0\}.$

(3) Assume that $L^*=[b,\infty)$ for some $b\in\mathbb{R}$ and
$U^*=\emptyset.$ We show that $U^*=(-\infty,a]$ for some
$-\infty<a\le b<\infty.$ By the symmetric version of statement (2)
of Proposition \ref{combinations of L star and U star lemma}, it is
immediate that $U^*\neq\{0\}.$

Since $L^*=[b,\infty),$ condition (a) or (c) of Proposition
\ref{lower bound theorem} must hold, with $\theta_4=\infty$. Thus,
$\Pi_{\xi,\eta}(A_3)=0,$ which implies that $\theta_3'=-\infty.$
Also, since $\theta_4=\infty,$ it must be that
$\Pi_{\xi,\eta}(A_4\setminus A_1)=0.$ Since $U^*\neq\emptyset,$ it
must be that $\Pi_{\xi,\eta}(A_1\cap A_4)=0,$ and so
$\Pi_{\xi,\eta}(A_4)=0.$ This implies that one of conditions (b) or
(c) of Proposition \ref{upper bound theorem} must hold, and so
$U^*=(-\infty,\theta_1']$ or $U^*=(-\infty,\theta_2']$ respectively.
Now, if condition (a) of Proposition \ref{lower bound theorem}
holds, then $L^*=[\theta_2,\infty).$ Note that Lemma
\ref{inequalities proposition} states that
$\theta_1'\le0\le\theta_2'\le\theta_2,$ and hence the result is
proved for either form of $U^*.$

Alternatively, if condition (c) of Proposition \ref{lower bound
theorem} holds, then $L^*=[\theta_1,\infty)$ where
$\theta_1>-\infty,$ which implies that $\Pi_{\xi,\eta}(A_1\setminus
A_4)\neq0.$ Hence, condition (b) of Proposition \ref{upper bound
theorem} must hold and $U^*=(-\infty,\theta_1'].$ Lemma
\ref{inequalities proposition} states that $\theta_1'\le\theta_1,$
and so we are done. \halmos
\end{pf}

\begin{pf}[Proposition \ref{equivalence proposition}]
\begin{enumerate}
\item[(1)$\Leftrightarrow$(2)]
Assume $L\cap U\neq\emptyset$ and let $z_1,z_2\in L\cap U.$ We show
$z_1=z_2\neq 0$. By Proposition \ref{lower bound theorem}, $z\in L$
iff $\eta-zW$ is increasing and by Proposition \ref{upper bound
theorem}, $z\in U$ iff $\eta-zW$ is decreasing. Thus,
$\eta-z_1W=\eta-z_2W=0,$ which implies $z_1W=z_2W.$ Since $\xi$ is
not zero, $W$ is not zero, and thus $z_1=z_2.$ Further, if
$z_1=z_2=0,$ then $\eta$ must be both increasing and decreasing,
which requires that $\eta$ be identically zero. Since we have
rejected this case, it must be that $z_1=z_2\neq 0.$

\item[(2)$\Leftrightarrow$(3)] Suppose $L\cap U=\{c\}.$ Then $V_t=c$ for all $t\ge 0$ whenever $V_0=c,$
which implies $e^{\xi_t}\left(c+Z_t\right)=c,$ which implies
$V_t=e^{\xi_t}(z-c)+c,$ as required. Conversely, suppose
$V_t=e^{\xi_t}(z-c)+c.$ Clearly, $c\in L\cap U$ and so $L\cap
U\neq\emptyset,$ which implies $L\cap U=\{c\}$ by the above.
\item[(2)$\Leftrightarrow$(4)]
By the definitions of $\delta$ and $\Upsilon,$ it is clear that $c$
is an absorbing point iff $\delta(c)=\Upsilon(c)=c,$ and the
definitions of $L$ and $U$ imply that this occurs iff $c\in L\cap
U.$

\item[(2)$\Rightarrow$(5)]
Assume $L\cap U=\{c\}$ where $c\neq 0.$ Propositions \ref{lower
bound theorem} and Proposition \ref{upper bound theorem} immediately
imply that equation (\ref{covariance matrix equation}) is satisfied
for $u=c,$ and imply respectively that $g(c)\ge 0$ and $g(c)\le 0$,
thus giving $g(c)=0.$ Finally, since $(2)\Rightarrow(3),$ the
equation $Z_t:=\int_0^t e^{-\xi_{s-}}\ud
\eta_s=c\left(e^{-\xi_t}-1\right)$ holds, which implies that
$e^{-\xi_{t-}}\Delta\eta_t=c\left(e^{-\xi_t}-1\right)-c\left(e^{-\xi_{t-}}-1\right)$
and so $\Delta\eta_t=c\left(e^{-\Delta\xi_t}-1\right).$

\item[(5)$\Rightarrow$(2)]
Assume that the conditions of statement (5) hold for $c\neq 0.$ We
prove $c\in L.$ Since (\ref{covariance matrix equation}) is
satisfied for $u=c,$ and $g(c)=0$ holds, we know that conditions (i)
and (iii) of Proposition \ref{lower bound theorem} are respectively
satisfied for $u=c.$ Thus it suffices to prove condition (ii) of
Proposition \ref{lower bound theorem} is satisfied for $u=c,$ or
equivalently, show $c\in L^*.$ If $\Pi_{\xi,\eta}=0$ then this is
trivial since $L^*=\mathbb{R}.$ Now suppose that $\Pi_{\xi,\eta}$ is
supported on the curve $\{(x,y):y-c(e^{-x}-1)=0\}$ for
$c\in\mathbb{R}.$ If $c> 0,$ $\Pi_{\xi,\eta}(A_2)\neq 0$ and
$\Pi_{\xi,\eta}(A_4)\neq 0,$ then $\theta_2=\theta_4=c$ and so
$L^*=\{c\}.$ If $c\ge 0,$ $\Pi_{\xi,\eta}(A_2)= 0$ and
$\Pi_{\xi,\eta}(A_4)\neq 0,$ then $\theta_2=0$ and $\theta_4=c,$ and
so $L^*=[0,c].$ If $c\ge 0,$ $\Pi_{\xi,\eta}(A_2)\neq 0$ and
$\Pi_{\xi,\eta}(A_4)= 0,$ then $\theta_2=c$ and $\theta_4=\infty,$
and so $L^*=[c,\infty).$ In each of these three cases, $c\in L^*.$
The proof for $c< 0$ is similar and we omit.

A symmetric argument proves that $c\in U.$ Hence, $c\in L\cap U$
which, by the equivalence of statements (1) and (2), implies that
$L\cap U=\{c\},$ as required.

\item[(2)$\Leftrightarrow$(6)]
$L\cap U=\{c\}$ iff $\eta-cW=0$ where $e^{-\xi_t}=\epsilon(W)_t$
which occurs iff $e^{-\xi_t}=\epsilon(\eta/c)_t$.
\end{enumerate}
Now assume that the above statements (1)-(6) hold. If
$\Sigma_{\xi,\eta}\neq0$ and both $L$ and $U$ are non-empty, then
Propositions \ref{lower bound theorem} and \ref{upper bound theorem}
immediately imply that $L=U=\{c\}$ where
$c=-\frac{\sigma_{\xi,\eta}}{\sigma^2_\xi}.$ For examples of \Levy
processes $(\xi,\eta)$ satisfying statements (1)-(6) and such that
$\xi$ drifts to $\infty$ a.s., $\xi$ drifts to $-\infty$ a.s. or
$\xi$ oscillates a.s., see Example \ref{second brownian example}.

If $\Sigma_{\xi,\eta}=0$ then the statements (a), (b) and (c) follow
immediately by examining the equation for $V$ in statement (3)
above. For examples of \Levy processes $(\xi,\eta)$ satisfying
statement (c) and such that $\xi$ drifts to $\infty$ a.s., $\xi$
drifts to $-\infty$ a.s. or $\xi$ oscillates a.s., see Example
\ref{degenerate example}.\halmos
\end{pf}

\begin{pf}[Theorem \ref{combinations proposition}] Assume that
$L\cap U=\emptyset.$ Suppose, firstly, that
$\Sigma_{\xi,\eta}\neq0.$ We must show that $(\xi,\eta)$ exists such
that (1), (2) or (3) occurs, and for each of these cases, we must
show that $\xi$ can satisfy each of the three asymptotic behaviours.
For case (1), this is obvious. Choosing $(\xi,\eta)$ such that
$\Sigma_{\xi,\eta}$ does not satisfy equation (\ref{covariance
matrix equation}) implies that $(\xi,\eta)$ fails both propositions,
and so $L=U=\emptyset,$ regardless of the choice of
$(\tilde{\gamma}_\xi,\tilde{\gamma}_\eta)$ and $\Pi_{\xi,\eta}.$
Clearly, we can make suitable choices for these objects to obtain
the desired asymptotic behaviour of $\xi.$ For case (2), our
existence claims are proven by Example \ref{first brownian example},
and case (3) is symmetric. It follows from Proposition
\ref{endpoints of L proposition}, and the symmetric version for $U,$
that whenever $L$ and $U$ are non-zero, they are each equal to
$\{-\sigma_{\xi,\eta}/\sigma^2_\xi\}.$ Hence, no cases, other than
(1), (2) and (3) of Theorem \ref{combinations proposition}, can
exist.

Now suppose that $\Sigma_{\xi,\eta}=0.$ We must show that
$(\xi,\eta)$ exists such that (a), (b) or (c) occurs, and for each
of these cases, we must show that $\xi$ can satisfy the specified
asymptotic behaviours. Examples \ref{first non-Brownian example} and
\ref{second non-Brownian example} present $(\xi,\eta)$ such that
$L=\emptyset,$ whilst $U$ may be of form $\emptyset,$ $\{a\}$ or
$[a,b]$ for $-\infty<a<b<\infty,$ and for each of these
combinations, it is shown that $\xi$ can satisfy the three
asymptotic behaviours. In Example \ref{third non-Brownian example},
$L=\emptyset,$ $U$ is of form $[b,\infty)$ for $b\in\mathbb{R},$ and
$\xi$ drifts to $-\infty$ a.s. In Example \ref{fifth non-Brownian
example}, $L=\emptyset,$ $U$ is of form $(-\infty,a]$ for
$a\in\mathbb{R},$ and $\xi$ drifts to $\infty$ a.s. These four
examples prove the existence claims for (a), and the case (b) is
symmetric. In Example \ref{fourth non-Brownian example},
$L=(-\infty,a],$ $U=[b,\infty)$ for $-\infty<a<b<\infty$ and $\xi$
drifts to $-\infty$ a.s. In Example \ref{sixth non-Brownian
example}, $U=(-\infty,a],$ $L=[b,\infty)$ for $-\infty<a<b<\infty,$
and $\xi$ drifts to $\infty$ a.s. These two examples prove the
existence claims for (c).

We now assume that $\Sigma_{\xi,\eta}=0,$ $L\neq\emptyset,$
$U\neq\emptyset$ and $L\cap U=\emptyset.$ We prove that no cases,
other than those listed in (c), can exist. As noted in point (2) of
Remark \ref{upper bound remarks}, it follows from our assumptions
that $(\xi,\eta)$ is finite variation and $g$ is linear.

Suppose that $L=[a,b]$ for some $-\infty<a\le b<\infty.$ We show
that this causes a contradiction with our assumptions. If
$L^*=[c,d]$ for some $-\infty<c\le a\le b\le d<\infty,$ then point
(2) of Proposition \ref{combinations of L star and U star lemma}
states that $U^*=\emptyset$ or $U^*=L^*=\{c\}=\{d\}.$ Thus,
$U=\emptyset$ or $U=L=\{a\}=\{b\},$ both of which contradict our
assumptions. Hence, it must be the case that $L^*=[c,\infty)$ for
some $-\infty<c\le a,$ or $L^*=(-\infty,d]$ for some $b\le
d<\infty.$

Thus, we suppose that $L=[a,b]$ and $L^*=[c,\infty)$ for some
$-\infty<c\le a\le b<\infty.$ The case in which $L^*=(-\infty,d]$
for some $b\le d<\infty$ is symmetric. We know $g(u)=d_\eta+ud_\xi.$
If $d_\xi\ge 0$ then it must be that $b=\infty,$ which we have
rejected. Hence $d_\xi<0,$ and we must have
$b=-\frac{d_\eta}{d_\xi}\ge a.$ Thus, since $U$ is non-empty, $L\cap
U=\emptyset,$ and $g(u)\le 0$ on $U,$ it must be that $U\subset
[b,\infty).$ However, point (3) of Proposition \ref{combinations of
L star and U star lemma} implies that $U^*\cap
[b,\infty)=\emptyset.$ Hence $U$ is empty, and we have a
contradiction. This completes the proof that $L\neq[a,b]$ for some
$-\infty<a\le b<\infty.$

We now assume that $L=[b,\infty)$ for $b\in\mathbb{R}.$ We first
prove that $\xi$ is a subordinator, which is another of the
statements of  Proposition \ref{endpoints of L proposition} and
point (2) of Remark \ref{lower bound remarks}, imply respectively,
that $(\xi,\eta)$ has no Brownian component, and $(\xi,\eta)$ is of
finite variation. Thus, we can write $g(u)=d_\eta+ud_\xi.$
Proposition \ref{lower bound theorem} implies that $g(u)\ge0$ on
$[b,\infty)$ and hence $d_\xi\ge 0.$ Finally, it must be that
$L^*=[c,\infty)$ for some $-\infty\le c\le b.$ It is a consequence
of the proofs of statements (1) and (3) of Proposition
\ref{combinations of L star and U star lemma}, that $\xi$ has no
negative jumps. Thus $\xi$ is a subordinator.

Now, we assume that $L=[b,\infty)$ for $b\in\mathbb{R}$ and
$U=\emptyset.$ We prove that $U=(-\infty,a]$ for some
$-\infty<a<b<\infty.$ Note that $L^*=[c,\infty)$ for some
$-\infty\le c\le b,$ so statement (3) of Proposition
\ref{combinations of L star and U star lemma} implies that
$U^*=(-\infty,d]$ for some $-\infty<d\le c.$ Since
$g(u)=d_\eta+ud_\xi$ and $d_\xi\ge 0,$ $U=(-\infty,a]$ for some
$-\infty<a\le d.$ Since we have assumed $L\cap U=\emptyset,$ $a<b$
as required.

If we assume that $U=(-\infty,a]$ for $a\in\mathbb{R},$ it can be
shown, using a method of proof similar to the one above, that $\xi$
is a subordinator, and $L=\emptyset$ or $L=[b,\infty)$ for some
$-\infty<a<b<\infty.$ We omit the details.

Now, if we assume $L=(-\infty,a]$ for $a\in\mathbb{R},$ then
symmetric proofs to the ones above, show that $-\xi$ is a
subordinator, and $U=\emptyset$ or $U=[b,\infty)$ for
$-\infty<a<b<\infty.$ Similarly, if we assume $U=[b,\infty)$ for
$b\in\mathbb{R},$ then symmetric proofs show that $-\xi$ is a
subordinator, and $L=\emptyset$ or $L=(-\infty,a]$ for
$-\infty<a<b<\infty.$ \halmos
\end{pf}

\begin{pf}[Proposition \ref{convergent divergent proposition}]
Assume $L\cap U=\emptyset.$ In the above proof of Theorem
\ref{combinations proposition}, it was shown that if $L=[b,\infty)$
for $b\in\mathbb{R}$ then $(\xi,\eta)$ is of finite variation,
$\Sigma_{\xi,\eta}=0,$ $d_\xi\ge 0,$ $\Pi_{\xi,\eta}(A_3)=0,$
$\Pi_{\xi,\eta}(A_4\setminus A_1)=0,$ and $\theta_2<\infty.$ It is
clear from Propositions \ref{lower bound theorem} and \ref{endpoints
of L proposition} that the converse also holds. A similar proof
shows that $U=(-\infty,a]$ for $a\in\mathbb{R}$ iff $(\xi,\eta)$ is
of finite variation, $\Sigma_{\xi,\eta}=0,$ $d_\xi\ge 0,$
$\Pi_{\xi,\eta}(A_4)=0,$ $\Pi_{\xi,\eta}(A_3\setminus A_2)=0,$ and
$\theta_1'>-\infty.$ Combining these two sets of iff conditions
immediately gives iff conditions for the case in which
$U=(-\infty,a]$ and $L=[b,\infty)$ with $-\infty<a<b<\infty.$ Since
$V$ is increasing on $L$ and decreasing on $U,$ and $V$ is a strong
Markov process, it is clear that in this situation
$\lim_{t\rightarrow\infty}|V_t|=\infty$ a.s. for any
$V_0=z\in\mathbb{R}.$

It follows by symmetric methods that $L=(-\infty,a]$ and
$U=[b,\infty)$ for $-\infty<a<b<\infty$ iff the stated conditions in
Proposition \ref{convergent divergent proposition} hold. The only
extra proof needed is to show that in this situation, $V$ is
strictly stationary. In \cite{LindnerMaller05} it is shown that
\[V_t=_De^{\xi_t}z+\int_0^te^{\xi_{s-}}\ud K^{\xi,\eta}_s.
\]
By Theorem 2 in \cite{EricksonMaller05} it is shown that if
$\lim_{t\rightarrow\infty}\xi_t=-\infty$ and and the integral
condition $I_{-\xi,K^{\xi,\eta}}=\infty$ holds, then
$|\int_0^te^{\xi_{s-}}\ud K^{\xi,\eta}_s|\rightarrow_P\infty$ as
$t\rightarrow\infty.$

As noted, if $L=(-\infty,a]$ and $U=[b,\infty)$ with
$-\infty<a<b<\infty$ then $-\xi$ is a subordinator and so
$\lim_{t\rightarrow\infty}\xi_t=-\infty$ a.s. Now if
$I_{-\xi,K^{\xi,\eta}}=\infty$ then by the above, and since
$\lim_{t\rightarrow\infty}e^{\xi_t}=-\infty$ a.s, it must be that
$|V_t|\rightarrow_D\infty.$ However this is impossible since $V$ is
increasing on $L$ and decreasing on $U.$ Thus, we must have
$I_{-\xi,K^{\xi,\eta}}<\infty.$ Hence, by Theorem 2.1 in
\cite{LindnerMaller05}, $V$ is strictly stationary and converges in
distribution to $\int_0^\infty e^{\xi_{s-}}\ud
K^{\xi,\eta}_s:=V_\infty.$ Since $V$ is increasing on $L$ and
decreasing on $U,$ and $V$ is a strong Markov process, it is clear
that $V_\infty$ has support $(a,b).$ \halmos
\end{pf}

\begin{pf}[Theorem \ref{linking theorem}] Assume $Z_t\rightarrow Z_\infty$
a.s. as $t\rightarrow\infty,$  where $Z_\infty$ is a finite random
variable. Suppose that for all $c\in\mathbb{R},$ equation
(\ref{degenerate case equation}) does not hold. This implies that
$Z_\infty$ is continuous.  As noted in Section \ref{intro section},
a necessary condition for the convergence of $Z_t,$ is
$\lim_{t\rightarrow\infty}\xi_t=\infty$ a.s., which implies that
$e^{\xi_t}\rightarrow\infty$ a.s. Since $Z_\infty$ is finite a.s.,
and $e^{\xi_t}\rightarrow\infty$ a.s., it is clear from the
definition $V_t:=e^{\xi_t}(z+Z_t),$ that
\begin{equation}\label{limit linking
equation}P(\lim_{t\rightarrow\infty}V_t=\infty|V_0=z)=P(Z_\infty>-z).
\end{equation}
Now let $a\le\sup U.$ By definition of $U,$
$P(\lim_{t\rightarrow\infty}V_t=\infty|V_0=a)=0$ which implies, by
equation (\ref{limit linking equation}), that $Z_\infty<-a$ a.s., as
required.

Conversely, let $a>\sup U.$ We prove $P(Z_\infty>-a)>0.$ Since we
have assumed that $|Z_\infty|<\infty$ a.s., we can choose $x>a$ such
that $P(Z_\infty>-x)>0.$ Note that $\Upsilon(a)=\infty$ and so there
exists a fixed time $T>0$ such that $P(V_T\ge x|V_0=a)>0.$

Hence, using (\ref{limit linking equation}), the law of conditional
probability and the Markov property,
\begin{eqnarray*}P(Z_\infty>-a)&=&P(\lim_{t\rightarrow\infty}V_t=\infty|V_0=a)\\
&\ge&P(\lim_{t\rightarrow\infty}V_t=\infty|V_T\ge x)P(V_T\ge
x|V_0=a)\\
&\ge&P(\lim_{t\rightarrow\infty}V_t=\infty|V_0=x)P(V_T\ge
x|V_0=a)\\
\end{eqnarray*}
which is greater than zero by (\ref{limit linking equation}) and the
choice of $x$ and $T.$ Thus,
\begin{equation}\label{repeat proposition statement}a\le\sup U~~\mathrm{iff}~~Z_\infty<-a~~\mathrm{a.s.}
\end{equation}
Now we prove $-\sup U=m$ where
$m:=\inf\{u\in\mathbb{R}|Z_\infty<u~\mathrm{a.s.}\}.$ By equation
(\ref{repeat proposition statement}), $Z_\infty<-\sup U$ and thus
$-\sup U\ge m.$ By assumption, $Z_\infty$ has no atoms and so
$Z_\infty<m$ a.s. Thus, equation (\ref{repeat proposition
statement}) implies that $-m\le\sup U.$ The proofs of the statements
for $L$ are symmetric.

Now assume that there exists $c\in\mathbb{R}$ such that equation
(\ref{degenerate case equation}) holds, and assume that
$Z_t\rightarrow Z_\infty$ a.s. as $t\rightarrow\infty.$ By equation
(\ref{degenerate case equation}) it is immediate that $Z_\infty=-c$
a.s. Further, since $\xi$ drifts to $\infty$ a.s., Proposition
\ref{equivalence proposition} implies that $L=U=\{c\},$ or
$U=(-\infty,c]$ and $L=[c,\infty).$ In both of these cases, $\inf
L=\sup U=c.$ \halmos
\end{pf}

\begin{pf}[Theorem \ref{minor certain ruin theorem}] (1) Assume $L\cap
U=\emptyset,$ $\sup U\ge 0$ and $L\cap[0,\sup U]=\emptyset,$ and let
$0\le u\le\sup U.$ We want to prove that $\psi(u)=1.$ Note that
there exists $z\ge u$ such that $z\in U,$ and so $\Upsilon(z)=z.$
Since $\psi(u)\ge\psi(z),$ it suffices to prove that $\psi(z)=1.$

Since $L\cap[0,\sup U]=\emptyset,$ we know $\delta(z)<0,$ which
implies that $P_z(\inf_{t>0}V_t<0)>0.$ Thus, there exists a fixed
time $T\in\mathbb{R}$ such that $P_z(\inf_{0<t\le T}V_t<0):=m>0.$
Let $n\in\mathbb{N}$ and let $A$ be the distribution of $V_{nT}$
conditional on both $V_0=z$ and $\inf_{0<t\le nT}V_t\ge0.$ Since
$\Upsilon(z)=z$ we know $A\le z$ a.s. Now
\[P_z\left(\inf_{nT<t\le
(n+1)T}V_t<0\Big|\inf_{0<t\le
nT}V_t\ge0\right)=P_A\left(\inf_{0<t\le T}V_t<0\right)\ge m,
\]
where the equality follows from the Markov property and the
inequality follows from the fact that $A\le z$ and $V_t$ is
increasing in $z.$ Define $P^n:=P_z\left(\inf_{0<t\le
nT}V_t<0\right)$ for all $n\in\mathbb{N}.$ By the law of total
probability
\[P^{n+1}=P^n+P_z\left(\inf_{nT<t\le (n+1)T}V_t<0\Big|\inf_{0<t\le
nT}V_t\ge0\right)\left(1-P^n\right)
\]
and so $P^{n+1}\ge P^n+(1-P^n)m$ where $P^1=m\in(0,1).$ This implies
that $P^n\ge 1-(1-m)^n$ which implies that
$\lim_{n\rightarrow\infty}P^n=1,$ and hence
$P_z\left(\inf_{0<t}V_t<0\right)=1$ by the continuity property of
measures.

(2) Assume $L\cap U=\emptyset,$ $\sup L\ge 0,$ and $U\cap[0,\sup
L]=\emptyset.$ We let $z\ge 0$ and prove that $\psi(z)<1.$ If
$z\ge\inf L$ then $\psi(z)=0$ by definition. Thus, it suffices to
assume $0\le z<\inf L.$

Suppose $\psi(z)=1.$ By assumption, $\Upsilon(z)>\inf L$ and so, by
definition, $P(C)>0$ where $C:=\{\sup_{t\ge 0} V_t\ge \inf L\}.$ By
definition of $L,$ $\lim_{t\rightarrow\infty}V_t\ge\inf L$ a.s. for
all $\omega\in C.$ Let $T_1:=\inf\{t>0|V_t<0\}$ and
$T_n:=\inf\{t>T_{n-1}|V_t<V_{T_{n-1}}\}$ for integers $n>1.$ By
assumption, $\psi(z)=1$ and so $T_1$ is finite a.s. Further, the
strong Markov property of $V$ implies that $\{T_n\}$ is a sequence
of stopping times increasing towards infinity as
$n\rightarrow\infty,$ and each $T_i$ is a.s. finite. In particular,
each $T_i$ is a.s. finite on $C.$ However $V_{T_n}<0$ a.s. which
contradicts the fact that $\lim_{t\rightarrow\infty}V_t>\inf L$ a.s.
on $C.$ Hence $\psi(z)<1.$ The proof of the case in which
$U\cap[0,\sup L]\neq\emptyset$ is almost identical, and we omit.
\halmos
\end{pf}

\begin{pf}[Theorem \ref{ruin prob theorem for L and U}](1): Assume $L\cap U=\emptyset,$ $\lim_{t\rightarrow\infty}\xi_t=-\infty$ a.s. and
$I_{-\xi,K^{\xi,\eta}}<\infty.$ Suppose that
$L\cap[0,\infty)\neq\emptyset.$ Since $\xi$ drifts to $-\infty$
a.s., Propositions \ref{equivalence proposition} and
\ref{combinations proposition} imply that one of conditions (a) or
(b) must hold. Further, it follows from statement (2) of Proposition
\ref{minor certain ruin theorem} and the definition of $L,$ that
$0<\psi(z)<1$ for all $0\le z<\inf L,$ and $\psi(z)=0$ for all
$z\ge\inf L.$

Now suppose that $L\cap[0,\infty)=\emptyset.$ We let $z\ge 0$ and
prove that $\psi(z)=1.$ Let $N$ be a Poisson process with parameter
$\lambda,$ let $D_i$ be an iid sequence of 1-dimensional exponential
random variables and let $C_i=1$ for all $i.$ Suppose that $N,$
$D_i$ and $(\xi,\eta)$ are mutually independent and define the
compound Poisson process $W_t:=\sum_{i=1}^{N_t}(C_i,D_i).$ Now
define a new \Levy process
$(\xi^\diamond_t,\eta^\diamond_t):=(\xi_t,\eta_t)+W_t,$ and denote
the associated GOU by $V^\diamond.$ For $V^\diamond,$ denote the
upper and lower bound functions, the sets of upper and lower bounds,
and the ruin probability function by $\Upsilon^\diamond,$
$\delta^\diamond,$ $U^\diamond$, $L^\diamond$ and $\psi^\diamond$
respectively.

Define $T_z:=\inf\{t>0:V_t<0|V_0=z\}.$ Since $\sup L<0,$ we know
$\delta(z)<0$ and hence $T_z$ is finite a.s. Note that
$V_0=V_0^\diamond=z.$ Also, whenever $V_{t-}\ge 0,$ every jump
$\Delta W_t$ causes a non-negative jump $\Delta V_t.$ Hence $V_t\le
V^\diamond_t$ a.s. on $t\le T_z.$ This implies that
$\psi(z)\ge\psi^\diamond(z).$ Thus it suffices to show that
$\psi^\diamond(z)=1$. To do this, we first need to prove that
$V^\diamond$ is strictly stationary.

We show that $\lambda>0$ can be chosen small enough such that
$\lim_{t\rightarrow\infty}\xi_t^\diamond=-\infty.$ Since
$\lim_{t\rightarrow\infty}\xi_t=-\infty$, either
$E(\xi_1)\in[-\infty,0)$ or $E(\xi_1)$ does not exist. If
$E(\xi_1)\in[-\infty,0)$ then $E(\xi_1^\diamond)=E(\xi_1)+\lambda$
and so we can choose $\lambda$ small enough such that
$E(\xi_1^\diamond)<0,$ which implies that
$\lim_{t\rightarrow\infty}\xi_t^\diamond=-\infty.$ If $E(\xi_1)$
does not exist then we know $E(\xi_1^\diamond)$ does not exist. We
show that $\lim_{t\rightarrow\infty}\xi_t^\diamond=-\infty$ holds
for any $\lambda>0.$ Note that $\xi^\diamond=\xi+N$ and , as noted
in Section \ref{intro section}, $J_\xi^+<\infty$ since $E(\xi_1)$
does not exist and $\lim_{t\rightarrow\infty}\xi_t=-\infty$.  Also
note that $\overline{\Pi}_{\xi^\diamond}^-=\overline{\Pi}_\xi^-$ and
so $A_{\xi^\diamond}^-=A_\xi^-.$ Since $\xi$ and $N$ are independent
we have $\overline{\Pi}_{\xi^\diamond}^+=\overline{\Pi}_\xi^+
+\overline{\Pi}_N^+.$ Further $\overline{\Pi}_N^+(x)=0$ for all
$x\ge 1.$ Hence $J_{\xi^\diamond}^+=J_\xi^+$ and so is finite. As
noted in Section \ref{intro section}, this implies that
$\lim_{t\rightarrow\infty}\xi_t^\diamond=-\infty.$

We now show that $(\xi^\diamond,\eta^\diamond)$ satisfies
$I_{-\xi^\diamond,K^{\xi^\diamond,\eta^\diamond}}<\infty.$ Since
$(\xi,\eta)$ and $W$ are independent, it is clear from the
definitions in Section \ref{intro section} that
$K^{\xi^\diamond,\eta^\diamond}_t=K^{\xi,\eta}_t+K^W_t$ and
$\overline{\Pi}_{K^{\xi^\diamond,\eta^\diamond}}(y)=\overline{\Pi}_{K^{\xi,\eta}}(y)+\overline{\Pi}_{K^W}(y).$
And, as above, $A_{-\xi^\diamond}^+=A_{-\xi}^+.$ Hence
\[I_{-\xi',K^{\xi^\diamond,\eta^\diamond}}=I_{-\xi,K^{\xi,\eta}}+\int_{(e,\infty)}\left(\frac{\ln(y)}{A_{-\xi}^+(\ln(y))}\right)|\overline{\Pi}_{K^{W}}(\ud
y)|.
\]
By the choice of $W$ it is clear that $K^W_1$ has a finite expected
value which implies that $\int_{(e,\infty)} y
|\overline{\Pi}_{K^{W}}(\ud y)|<\infty.$ Hence
$I_{-\xi',K^{\xi^\diamond,\eta^\diamond}}<\infty.$ Thus $V^\diamond$
is strictly stationary.

For a Lebesgue set $\Lambda$ define
$T_\Lambda^\diamond:=\inf\{t>0:V_t^\diamond\in\Lambda\}.$ Note that
$\theta_1'^\diamond=-\infty$ and hence Proposition \ref{upper bound
theorem} implies that $\Upsilon^\diamond(u)=\infty$ for all
$u\in\mathbb{R},$ or equivalently, $U^\diamond=\emptyset.$ Also,
$\theta_1^\diamond=0,$ and so Proposition \ref{lower bound theorem}
implies that $L^\diamond\cap(-\infty,0)=\emptyset,$ whilst the fact
that $L\cap(0,\infty)=\emptyset$ clearly implies that
$L'\cap(0,\infty)=\emptyset.$

These facts imply that, for all $a$ and $u$ in $\mathbb{R},$
$P\left(T_{(-\infty,a]}^\diamond<\infty|V_0^\diamond=u\right)>0$ and
$P\left(T_{[a,\infty]}^\diamond<\infty|V_0^\diamond=u\right)>0.$
Since $D$ is an exponential random variable, it is clear that
$V_t^\diamond$ has a continuous density with respect to Lebesgue
measure. Hence $P\left(T_\Lambda^\diamond<\infty\right)>0$ for any
set $\Lambda$ with positive Lebesgue measure. This result, and the
fact that $V^\diamond$ is strictly stationary, allows us to mimic
the argument of Theorem 3.1 (a) in Paulsen \cite{Paulsen98}. Let $S$
be an independent standard exponential variable and define the
resolvent kernel
\[K(z,\Lambda):=\int_0^\infty P_z(V^\diamond_t\in \Lambda)e^{-t}\ud t=P_z(V^\diamond_S\in
\Lambda).
\]
Proposition 2.1 of \cite{MeynTweedieII93} implies that $V^\diamond$
is $\phi$-irreducible for the measure $\phi=\lambda K.$ Using the
language of \cite{MeynTweedieII93} p.495 and 496, it is clear that
$K$ has a continuous nontrivial component for all $z$ and hence is a
T-process. Since $V^\diamond$ is strictly stationary it is clear
that $V^\diamond$ is non-evanescent, as defined in
\cite{MeynTweedieII93} p.494. Thus Theorem 3.2 of
\cite{MeynTweedieII93} p.494 implies that $V^\diamond$ is Harris
recurrent, as defined in \cite{MeynTweedieII93} p490, which clearly
implies that $\psi^\diamond(z)=1$ as required.

(2) Assume that $L\cap U=\emptyset,$ $E(\xi_1)=0,$
$E(e^{|\xi_1|})<\infty$ and there exist $p,q>1$ with $1/p+1/q=1$
such that $E\left(e^{-p\xi_1}\right)<\infty$ and
$E\left(|\eta_1|^q\right)<\infty.$

Suppose that $L\cap[0,\infty)\neq\emptyset.$ Since $\xi$ oscillates
a.s., Proposition \ref{combinations proposition} implies that
$L=[a,b]$ and $U=\emptyset$ where $-\infty<a\le b<\infty$ and $b\ge
0.$ Hence, it follows from statement (2) of Proposition \ref{minor
certain ruin theorem} and the definition of $L,$ that $0<\psi(z)<1$
for all $0< z<a$ and $\psi(z)=0$ for all $z\ge a.$

Now suppose that $L\cap[0,\infty)=\emptyset.$ We let $z\ge 0$ and
prove that $\psi(z)=1.$ We know that
$P\left(\inf_{t>0}V_t<0|V_0=z\right)>0.$ However, it is possible
that for some $z>0,$ $P(V_1<0|V_0=z)=0.$ For example, this would
happen if $(\xi,\eta)$ has no Brownian component and $\sup L^*>0.$
Let $0=T_0<T_1<T_2<\ldots$ be random times such that $T_i-T_{i-1}$
are iid with exponential distribution and parameter $\lambda.$ Since
$T_1$ has infinite support it is clear that $\sup L<0$ implies
$P\left(V_{T_1}<0|V_0=z\right)>0$ for all $z\ge 0.$ Equation
(\ref{GOU definition}) implies that a.s.
\[V_{T_n}=e^{\xi_{T_n}-\xi_{T_{n-1}}}\left(e^{\xi_{T_{n-1}}}\left(z+\int_0^{T_{n-1}}e^{-\xi_{s-}}\ud\eta_s\right)\right)
+e^{\xi_{T_n}}\int_{T_{n-1}}^{T_n}e^{-\xi_{s-}}\ud\eta_s.\] Thus, if
we define $A_n:=e^{\xi_{T_n}-\xi_{T_{n-1}}},$
$B_n:=e^{\xi_{T_n}}\int_{T_{n-1}}^{T_n}e^{-\xi_{s-}}\ud\eta_s$ and
the stochastic difference equation $W_n:=A_nW_{n-1}+B_n$ with
$W_0:=V_0=z$ then $W_n=V_{T_n}$ a.s. for all $n\in\mathbb{N}.$ Note
that the term $e^{\xi_{T_n}}$ in $B_n$ cannot be brought under the
integral sign because it is not predictable. Since a \Levy process
has independent increments it is clear that $(A_n,B_n)$ is an
independent sequence. Now,
\begin{eqnarray*}(A_2,B_2)&=&\left(e^{\xi_{T_2}-\xi_{T_{1}}},e^{\xi_{T_2}-\xi_{T_{1}}}e^{\xi_{T_1}}\int_{T_1}^{T_2}e^{-\xi_{s-}}\ud\eta_s\right)\\
&=&\left(e^{\xi_{T_2}-\xi_{T_{1}}},e^{\xi_{T_2}-\xi_{T_{1}}}\int_{T_1}^{T_2}e^{-\left(\xi_{s-}-\xi_{T_1}\right)}\ud\eta_s\right)\\
&=&\left(e^{\xi_{T_2}-\xi_{T_{1}}},e^{\xi_{T_2}-\xi_{T_{1}}}\int_{T_1}^{T_2}e^{-\left(\xi_{s-}-\xi_{T_1}\right)}\ud(\eta_s-\eta_{T_1})\right)\\
&=&_D\left(e^{\xi_{T_1}},e^{\xi_{T_1}}\int_{T_1}^{T_2}e^{-\xi_{s-{T_1}}}\ud\eta_{s-T_1}\right)\\
&=&\left(e^{\xi_{T_1}},e^{\xi_{T_1}}\int_0^{T_1}e^{-\xi_{s-}}\ud\eta_s\right)=(A_1,B_1),
\end{eqnarray*}
where the second equality holds because $e^{\xi_{T_1}}$ is
predictable with respect to the integral, the third equality holds
because a \Levy process has identically distributed increments and
the final equality is obtained using a change of variables. The
argument for general $n$ is identical, and thus $(A_n,B_n)$ is an
iid sequence.

Now Proposition 1.1 and Corollary 4.2 of
\cite{BabillotBougerolElie97} state that if $P(A_1z+B_1=z)<1$ for
all $z\in\mathbb{R},$ $E(\ln A_1)=0,$ $A_1\not\equiv 1$ and there
exists $\alpha>0$ such that
\begin{equation}\label{Babillot condition}E\left(\left(|\ln
A_1|+\ln^+|B_1|\right)^{2+\alpha}\right)<\infty\end{equation} then
the discrete stochastic process $W$ has an invariant unbounded Radon
measure $\mu$ unique up to a constant factor such that the sample
paths $W_n,$ with $W_0=z,$ visit every open set of positive
$\mu$-measure infinitely often with probability 1, for every
$z\in\mathbb{R}.$ The first of these conditions follows from our
assumption that $L\cap U=\emptyset,$ using Proposition
\ref{equivalence proposition}. The second and third conditions
follow respectively from our assumptions that $E(\xi_1)=0,$ and
$\xi_1$ is not identically zero. We will
show later that our moment conditions on $\xi$ and $\eta$ ensure
equation (\ref{Babillot condition}) holds. Note that the Babillot
result implies that $\psi(z)=1$ if we can show
$\mu\left((-\infty,0)\right)>0.$ However by the definition of an
invariant measure,
\[\mu\left((-\infty,0)\right)=\int_{z\in\mathbb{R}}P(A_1z+B_1<0)\mu(\ud z)\ge\int_{z\in\mathbb{R}}P(V_{T_1}<0|V_0=z)\mu(\ud z).
\]
Thus if $\mu\left([0,\infty)\right)>0$ then
$\mu\left((-\infty,0)\right)>0$ since
$P\left(V_{T_1}<0|V_0=z\right)>0$ for all $z\ge 0.$ And if
$\mu\left([0,\infty)\right)=0$ then $\mu\left((-\infty,0)\right)>0$
since $\mu(\mathbb{R})>0.$ Thus we are done if we can prove equation
(\ref{Babillot condition}).

To do this, it suffices to assume $T_1=1$ and
$(A_1,B_1):=\left(e^{\xi_1},e^{\xi_1}\int_0^1
e^{-\xi_{s-}}\ud\eta_s\right)$ since we can choose the parameter
$\lambda$ of the increments to be arbitrarily small.  Note that if
$x,y>0$ and $\alpha>0$ then there exists $c_1>0$ such that
$(x+y)^\alpha\le c_1\left(x^\alpha+y^\alpha\right).$ Also
$\ln^+(x+y)\le\ln^+(x)+\ln^+(y)$ and
$\ln^+(xy)\le\ln^+(x)+\ln^+(y).$ Finally note that whenever
$0<\alpha\le 1$ there exists $c_2>0$ such that
$\ln^+(x)^{2+\alpha}\le c_2 x^\alpha.$ Using these four inequalities
it is clear that equation (\ref{Babillot condition}) is satisfied
whenever there exists $0<\alpha\le 1$ such that
$E\left(e^{\alpha\xi_1}\right)<\infty,$
$E\left(|\xi_1|^{2+\alpha}\right)<\infty$ and $E\left(\left|\int_0^1
e^{-\xi_{s-}}\ud\eta_s\right|^\alpha\right)<\infty.$ By Proposition
\ref{sup proposition}, and the fact that the existence of an
absolute exponential moment implies the existence of absolute
moments of all orders, the assumed moment conditions imply that
these conditions are satisfied for $\alpha=1.$

(3) Assume that $\lim_{t\rightarrow\infty}\xi_t=\infty$ a.s. and
$I_{\xi,\eta}<\infty.$ Suppose that $-\infty\le\sup U<z.$ Assume,
for the sake of contradiction, that $\psi(z)=1.$ Theorem
\ref{linking theorem} implies that $P(C)>0$ where
$C:=\{Z_\infty>-z\}.$ Since $\lim_{t\rightarrow\infty}\xi_t=\infty$,
we know that $\lim_{t\rightarrow\infty}V_t=\infty$ a.s. on $C.$ Now,
the same strong Markov property argument used in the proof of
statement (2) of Theorem \ref{minor certain ruin theorem}, gives a
contradiction. Hence $\psi(z)<1.$

Now suppose $U\cap[0,\infty)\neq\emptyset.$ Since $\xi$ drifts to
$\infty$ a.s., Theorem \ref{combinations proposition} implies that
either $U=[a,b]$ and $L=\emptyset$ where $-\infty\le z\le b<\infty$
and $b\ge 0,$ or $U=(-\infty,a]$ and $L=[b,\infty)$ for some $0\le
a<b<\infty.$ In both of these cases, statement (1) of Theorem
\ref{minor certain ruin theorem} implies that $\psi(z)=1$ for all
$z\le\sup U.$ Using the definition of $L,$ and the above result, it
is clear that $0<\psi(z)<1$ for all $\sup U<z<\inf L$ and
$\psi(z)=0$ for all $z\ge\sup L.$  \halmos
\end{pf}

\begin{pf}[Proposition \ref{degenerate certain ruin prop}] Assume
that $V_t=e^{\xi_t}(z-c)+c.$ By definition of $L,$ if $c\ge 0$ then
$\psi(z)=0$ for all $z\ge c.$

Let $0\le z<c.$ If $\xi$ drifts to $-\infty$ a.s. then
$\lim_{t\rightarrow\infty}V_t=c$ a.s. Thus, the strong Markov
property of $V$ implies that $\psi(z)<1,$ using a proof similar to
that used for statement (2) of Theorem \ref{minor certain ruin
theorem}. If $\xi$ oscillates a.s. then
$-\infty=\liminf_{t\rightarrow\infty}V_t<\limsup_{t\rightarrow\infty}V_t=c,$
and so $\psi(z)=1.$ If $\xi$ drifts to $\infty$ a.s. then
$\lim_{t\rightarrow\infty}V_t=-\infty$ a.s. which implies
$\psi(z)=1.$

Let $c<0\le z$. If $\xi$ drifts to $-\infty$ a.s. then
$\lim_{t\rightarrow\infty}V_t=c$ a.s. and so $\psi(z)=1.$ If $\xi$
oscillates a.s. then
$c=\liminf_{t\rightarrow\infty}V_t<\limsup_{t\rightarrow\infty}V_t=\infty,$
and so $\psi(z)=1.$ If $\xi$ drifts to $\infty$ a.s. then
$\lim_{t\rightarrow\infty}V_t=\infty$ a.s. which implies
$\psi(z)<1,$ using a strong Markov property argument. \halmos
\end{pf}

\begin{pf}[Theorem \ref{paulsen theorem 1}] Suppose that for all $c\in\mathbb{R}$ the degenerate case (\ref{degenerate case equation}) does not hold.
Then, by Proposition \ref{equivalence proposition}, $L\cap
U=\emptyset.$ It follows immediately from Theorem \ref{ruin prob
theorem for L and U} that $0<\psi(z)<1$ iff $0\le z<m<\infty$
whenever the assumptions for statement (1), or statement (2), of
Theorem \ref{paulsen theorem 1} are satisfied. Now suppose that
there exists $c\in\mathbb{R}$ such that equation (\ref{degenerate
case equation}) holds. Then it follows immediately from Proposition
\ref{degenerate certain ruin prop} that $0<\psi(z)<1$ iff $0\le
z<m<\infty$ whenever the assumptions for statement (1), or statement
(2), of Theorem \ref{paulsen theorem 1} are satisfied. In both these
situations, $m=c.$ \halmos
\end{pf}

\begin{pf}[Theorem \ref{paulsen theorem 2}] Assume $\lim_{t\rightarrow\infty}\xi_t=\infty$ a.s. and
$I_{\xi,\eta}<\infty.$ Assume that for all $c\in\mathbb{R}$ equation
(\ref{degenerate case equation}) does not hold, or equivalently,
$L\cap U=\emptyset.$ Theorem \ref{paulsen theorem 2} claims that
$\psi(0)=1$ iff $-\eta$ is a subordinator, or there exists $z>0$
such that $\psi(z)=1.$ This claim follows by combining two known
results: $\psi(z)=1$ iff $\sup U\ge 0$ and $z<\sup U,$ which is
implied by statement (3) of Theorem \ref{minor certain ruin
theorem}; secondly, $0\in U$ iff $-\eta$ is a subordinator, which is
stated in Proposition \ref{upper bound theorem}.

Theorem \ref{paulsen theorem 2} also states conditions on the
characteristic triplet of $(\xi,\eta)$ and claims these are
equivalent to the fact that there exists $z>0$ such that
$\psi(z)=1.$ However, using statement (3) of Theorem \ref{minor
certain ruin theorem}, we know there exists $z>0$ such that
$\psi(z)=1$ iff $\sup U>0.$ And Proposition \ref{upper bound
theorem} gives iff conditions on the characteristic triplet of
$(\xi,\eta)$ for the case $\sup U>0$ to occur. These conditions are
precisely the conditions stated in Theorem \ref{paulsen theorem 2}.

Finally, statements (1) and (2) of Theorem \ref{paulsen theorem 2}
contain values for $\sup\{z\ge 0:\psi(z)=1\}.$ However, these follow
from the unstated parallel version of Proposition \ref{endpoints of
L proposition} which gives exact values for the endpoints of $U.$

Now, assume that there exists $c\in\mathbb{R}$ such that the
degenerate equation (\ref{degenerate case equation}) holds, and
$L=U=\{c\}$. Since $\xi$ drifts to $\infty$ a.s., Proposition
\ref{equivalence proposition} implies that $\sup U=c.$ Thus,
Proposition \ref{degenerate certain ruin prop} implies that
$\psi(z)=1$ iff $\sup U\ge 0$ and $z<\sup U.$ Theorem \ref{paulsen
theorem 2} is proved for the degenerate case by combining this
statement with Proposition \ref{upper bound theorem} Proposition
\ref{upper bound theorem} and the parallel version of Proposition
\ref{endpoints of L proposition}, in an identical manner to the
above. The only difference is that the set $\{z\ge 0:\psi(z)=1\}$
does not contain its supremum in the degenerate case, since
$\sup\{z\ge 0:\psi(z)=1\}=U=L,$ and is an absorbing point. \halmos
\end{pf}

\subsection{Examples}
Propositions \ref{equivalence proposition}, \ref{combinations
proposition} and \ref{convergent divergent proposition} claim that
\Levy processes $(\xi,\eta)$ exist which satisfy particular
combinations of $L$ and $U,$ and particular asymptotic behaviour for
$\xi$. We now present examples which prove these claims. We use the
simplest \Levy processes possible. The \Levy measures will always be
finite activity, namely $\Pi_{\xi,\eta}(\mathbb{R}^2)<\infty.$
Hence, we can write $(\xi,\eta)$ in the form
$(\xi,\eta)_t=(d_\xi,d_\eta)t+(B_{\xi,t},B_{\eta,t})+\sum_{i=1}^{N_t}Y_i$
where $(B_{\xi,t},B_{\eta,t})$ is Brownian motion with covariance
matrix $\Sigma_{\xi,\eta},$ $N$ is a Poisson process with parameter
$\Lambda$ and $\{Y_i\}_{i=1}^\infty$ is an iid sequence of two
dimensional random variables with distribution $Y.$

\textbf{Examples with Brownian component} The first example is of a
\Levy process $(\xi,\eta)$ for which $L=\{a\},$ $U=\emptyset.$ The
second example is of a \Levy process for which $L=U=\{a\}.$ For both
examples we show how variables can be chosen so that $\xi$ drifts to
$\infty$ a.s., $\xi$ drifts to $-\infty$ a.s. or $\xi$ oscillates
a.s.

\begin{exmp}\label{first brownian example}Let
$(\xi,\eta)_t:=(d_\xi,2)t+(B_t,B_t)+\sum_{i=1}^{N_t}Y_i$ where $B$ is
a one-dimensional Brownian motion with variance $1,$ and
$P(Y=\left(10,10)\right)=1/2$ and $P(Y=\left(-10,10)\right)=1/2.$
The covariance matrix equation (\ref{covariance matrix equation}) is
satisfied for $u=-1.$ Condition (ii) of Proposition \ref{lower bound
theorem} is satisfied for $u=-1$, whilst condition (ii) of
Proposition \ref{upper bound theorem} is not satisfied. By equation
(\ref{linear g equation}), $g(-1)=3/2-d_\xi ,$ and so choosing $d_\xi
\le3/2$ implies that $L=-1$ and $U=\emptyset.$ However
$E(\xi_1)=d_\xi $ so if $0<d_\xi <3/2$ then $\xi$ drifts to $\infty$
a.s., if $d_\xi <0$ then $\xi$ drifts to $-\infty$ a.s., and if
$d_\xi =0$ then $\xi$ oscillates a.s.
\end{exmp}

\begin{exmp}\label{second brownian example}Let $(\xi,\eta)_t:=(d_\xi ,d_\eta )t+(B_t,-B_t).$
Equation (\ref{covariance matrix equation}) is satisfied for $u=1,$
whilst condition (ii) of Proposition \ref{lower bound theorem} and
condition (ii) of Proposition \ref{upper bound theorem} are
satisfied trivially. Equation (\ref{linear g equation}) implies
$g(1)=d_\eta +d_\xi -1/2.$ Thus, choosing $d_\xi =1/2-d_\eta $
implies that $L=U=1.$ Note $E(\xi_1)=d_\xi ,$ so if $d_\eta <1/2$
then $\xi$ drifts to $\infty$ a.s., if $d_\eta
>1/2$ then $\xi$ drifts to $-\infty$ a.s., and if $d_\eta =1/2$ then
$\xi$ oscillates a.s.
\end{exmp}

\textbf{Examples with no Brownian component} We now present seven
examples of \Levy processes $(\xi,\eta)$ with no Brownian component.
In Example \ref{degenerate example}, $L=U=\{a\}$ and we indicate how
the parameters can be changed in order to obtain each of the three
asymptotic behaviours for $\xi$. In Examples \ref{first non-Brownian
example} and \ref{second non-Brownian example}, $L=\emptyset,$
whilst $U$ may be of form $\emptyset,$ $\{a\}$ or $[a,b]$ for
$-\infty<a<b<\infty.$ We indicate how parameters can be changed in
order to obtain these different sets, and for each set, to obtain
the three possible asymptotic behaviours for $\xi.$ In Example
\ref{third non-Brownian example}, $L=\emptyset$ whilst $U$ is of
form $[b,\infty)$ for $b\in\mathbb{R}.$ In Example \ref{fourth
non-Brownian example}, $L=(-\infty,a]$ and $U=[b,\infty)$ for
$-\infty<a<b<\infty.$ For both these examples we show that $\xi$
drifts to $-\infty$ a.s. In Example \ref{fifth non-Brownian
example}, $L=\emptyset$ whilst $U$ is of form $(-\infty,a]$ for
$a\in\mathbb{R}.$ In Example \ref{sixth non-Brownian example},
$U=(-\infty,a]$ and $L=[b,\infty)$ for $-\infty<a<b<\infty.$ For
both these examples we show that $\xi$ drifts to $\infty$ a.s.

\begin{exmp}\label{degenerate example}
Let $(\xi,\eta)_t:=(d_\xi ,d_\eta )t+\sum_{i=1}^{N_t}Y_i$ where
$P(Y=(3,2e^{-3}-2))=1/2$ and $P(Y=(-3,2e^{3}-2))=1/2.$ Then
$\theta_2=\theta_2'=\theta_4=\theta_4'=2$ and $L^*=U^*=\{2\}.$ Note
that $g(u)=d_\eta +ud_\xi,$ so choosing $d_\eta =-2d_\xi $ implies
that $g(2)=0$ and hence $L=U=\{2\}.$ Since $E(\xi_1)=d_\xi ,$
choosing $d_\xi >0,$ $d_\xi <0,$ and $d_\xi =0,$ implies that $\xi$
drifts to $\infty$ a.s., $\xi$ drifts to $-\infty$ a.s. and $\xi$
oscillates a.s., respectively.
\end{exmp}

\begin{exmp}\label{first non-Brownian example}
Let $(\xi,\eta)_t:=(d_\xi ,d_\eta )t+\sum_{i=1}^{N_t}Y_i$ where
$P(Y=\left(4,-2)\right)=1/3$ and $P(Y=\left(-2,-3)\right)=1/3$ and
$P(Y=\left(-2,1)\right)=1/3.$ Then $L=\emptyset$ since
$\Pi_{\xi,\eta}(A_2)$ and $\Pi_{\xi,\eta}(A_3)$ are both non-zero,
whilst
$U^*=[\theta_4',\theta_2']=[\frac{-2}{e^{-4}-1},\frac{1}{e^2-1}]\cong[0.2,2].$
Now $U=\{u\in U^*:g(u)\le0\}$ and $g$ simplifies to $g(u)=d_\eta
+ud_\xi $. Note that $E(\xi_1)=d_\xi .$

Choosing $d_\xi =0$ and $d_\eta >0$ implies that $U=\emptyset$ and
$\xi$ oscillates a.s. Choosing $d_\xi >0$ and $d_\eta
>-\theta_4'd_\xi $ implies that $U=\emptyset$ and $\xi$ drifts to
$\infty$ a.s. Choosing $d_\xi <0$ and $d_\eta >-\theta_2'd_\xi $
implies that $U=\emptyset$ and $\xi$ drifts to $-\infty$ a.s.

Choosing $d_\xi =0$ and $d_\eta <0$ implies that $U=U^*\cong[0.2,2]$
and $\xi$ oscillates a.s.  Choosing $d_\xi >0$ and $d_\eta
<-\theta_2'd_\xi $ implies that $U=U^*\cong[0.2,2]$ and $\xi$ drifts
to $\infty$ a.s. Choosing $d_\xi <0$ and $d_\eta <-\theta_4'd_\xi $
implies that $U=U^*\cong[0.2,2]$ and $\xi$ drifts to $-\infty$ a.s.

Choosing $d_\xi >0$ and $d_\eta =-\theta_4'd_\xi $ implies that
$U=\{\theta_4'\}\cong\{0.2\}$ and $\xi$ drifts to $\infty$ a.s.
Choosing $d_\xi <0$ and $d_\eta =-\theta_2'd_\xi $ implies that
$U=\{\theta_2'\}\cong\{2\}$ and $\xi$ drifts to $-\infty$ a.s.
\end{exmp}
Note that for Example \ref{fifth non-Brownian example}, no
adjustment of $d_\xi $ and $d_\eta $ can result in $U=\{a\}$ with
$\xi$ oscillating a.s. We now present a different example with this
behaviour.

\begin{exmp}\label{second non-Brownian example}Let
$(\xi,\eta)_t:=(0,-2)t+\sum_{i=1}^{N_t}Y_i$ where
$P(Y=\left(2,e^{-2}-1)\right)=1/3$ and
$P(Y=\left(-1,e-1)\right)=1/3$ and $P(Y=\left(-1,-2)\right)=1/3.$
Then $L=\emptyset,$ $\theta_2=\theta_2'=\theta_4=\theta_4'=1,$ and
$U^*=\{1\}.$ Since $g$ simplifies to $g(u)=-2$ for all
$u\in\mathbb{R}$ we obtain $U=\{1\}.$ Since $E(\xi_1)=0,$ $\xi$
oscillates a.s.
\end{exmp}

\begin{exmp}\label{third non-Brownian example}Let
$(\xi,\eta)_t:=(0,-2)t+\sum_{i=1}^{N_t}Y_i$ where $P(Y=(-1,2))=1/3$
and $P(Y=(-2,-3))=1/3$ and $P(Y=(0,-5))=1/3.$ Then $L^*=\emptyset$
whilst
$U^*=[\theta_4',\theta_2']=[\frac{2}{e-1},\infty)\cong[1.2,\infty).$
Since $g(u)=-2$ for all $u\in\mathbb{R}$ we obtain $L=\emptyset$ and
$U=U^*$ Since $E(\xi_1)=-1.5,$ $\xi$ drifts to $-\infty$ a.s.
\end{exmp}

\begin{exmp}\label{fourth non-Brownian example}Let
$(\xi,\eta)_t:=(d_\xi ,d_\eta )t+\sum_{i=1}^{N_t}Y_i$ where
$P(Y=(-1,2))=1/2$ and $P(Y=(-2,-3))=1/2.$ Then
$L^*=[\theta_1,\theta_3]=(-\infty,\frac{-3}{e^2-1}]\cong(-\infty,-0.5]$
and
$U^*=[\theta_4',\theta_2']=[\frac{2}{e-1},\infty)\cong[1.2,\infty).$
Note that $g$ simplifies to $g(u)=d_\eta +ud_\xi $ and hence
choosing $d_\xi \le0$ and $d_\eta =0$ gives $L=L^*$ and $U=U^*.$
Since $E(\xi_1)=-1.5+d_\xi ,$ $\xi$ drifts to $-\infty$ a.s.
\end{exmp}

\begin{exmp}\label{fifth non-Brownian example}Let
$(\xi,\eta)_t:=\sum_{i=1}^{N_t}Y_i$ where $P(Y=(1,2))=1/3$ and
$P(Y=(1,8))=1/3$ and $P(Y=(0,-5))=1/3.$ Then $L^*=\emptyset$ whilst
$U^*=[\theta_3',\theta_1']=(-\infty,\frac{8}{e^{-1}-1}]\cong(-\infty,-12.6].$
Note that $g(u)=0$ for all $u\in\mathbb{R}$ so $L=L^*$ and $U=U^*.$
Since $E(\xi_1)=1,$ $\xi$ drifts to $\infty$ a.s.
\end{exmp}

\begin{exmp}\label{sixth non-Brownian example}Let
$(\xi,\eta)_t:=\sum_{i=1}^{N_t}Y_i$ where $P(Y=(1,2))=1/2$ and
$P(Y=(1,8))=1/2.$ Then
$L^*=[\theta_1,\theta_4]=[\frac{2}{e^{-1}-1},\infty)\cong[-3.2,\infty)$
and
$U^*=[\theta_3',\theta_1']=(-\infty,\frac{8}{e^{-1}-1}]\cong(-\infty,-12.6].$
Note that $g(u)=0$ for all $u\in\mathbb{R}$ so $L=L^*$ and $U=U^*.$
Since $E(\xi_1)=1,$ $\xi$ drifts to $\infty$ a.s.
\end{exmp}

\textbf{Acknowledgements.} We are very grateful to Professor Ross
Maller.

\bibliography{bibliography}
\bibliographystyle{elsart-num-sort}

\end{document}